\documentclass[12pt]{article}
\usepackage{a4wide}
\usepackage{amsmath}
\usepackage{amssymb}
\usepackage{booktabs}
\usepackage{comment}
\usepackage[textwidth=30mm]{todonotes}
\usepackage{latexsym}
\usepackage{mathpazo}
\usepackage{mathrsfs}
\usepackage{tikz}
\usepackage{booktabs}
\usepackage{hyperref}

\newtheorem{thm}{Theorem}[section]

\newtheorem{cor}[thm]{Corollary}

\newtheorem{exl}[thm]{Example}
\newtheorem{rem}[thm]{Remark}
\newtheorem{lemma}[thm]{Lemma}
\newtheorem{prop}[thm]{Proposition}

\numberwithin{equation}{section}

\newcommand{\N}{\mathbb{N}}

\newcommand{\A}{\mathcal{A}}

\newcommand{\F}{\mathscr{F}}
\newcommand{\I}{\mathscr{I}}
\newcommand{\J}{\mathcal{J}}
\newcommand{\K}{\mathcal{K}}

\newenvironment{proof}{\vspace{-2pt}\noindent\textbf{Proof:}}{$\Box$\\\medskip}
\newcounter{step}
\newenvironment{step}[1]{\refstepcounter{step}\noindent\textsc{Step \thestep{\rm #1}:}}{$\Box$\\}

\author{J\'{a}nos Flesch\footnote{Maastricht University, email: \texttt{j.flesch@maastrichtuniversity.nl}}\ \and Arkadi Predtetchinski\footnote{Maastricht University, email: \texttt{a.predtetchinski@maastrichtuniversity.nl}}\ \and William D Sudderth\footnote{University of Minnesota, email: \texttt{bill@stat.umn.edu}}\ \and Xavier Venel\footnote{LUISS, email: \texttt{xvenel@luiss.it}}}

\begin{document}
\title{Zero-one Laws for a Control Problem with Random Action Sets}
\maketitle
\begin{abstract}
\noindent\textbf{Abstract:} In many control problems there is only limited information about the actions that will be available at future stages. We introduce a framework where the Controller chooses actions $a_{0}, a_{1}, \ldots$, one at a time. Her goal is to maximize the probability that the infinite sequence $(a_{0}, a_{1}, \ldots)$ is an element of a given subset $G$ of $\N^\N$. The set $G$, called the goal, is assumed to be a Borel tail set. The Controller's choices are restricted: having taken a sequence $h_{t} = (a_{0}, \ldots, a_{t-1})$ of actions prior to stage $t \in \N$, she must choose an action $a_{t}$ at stage $t$ from a non-empty, finite subset $A(h_{t})$ of $\N$. The set $A(h_{t})$ is chosen from a distribution $p_{t}$, independently over all $t \in \N$ and all $h_{t} \in \N^{t}$. We consider several information structures defined by how far ahead into the future the Controller knows what actions will be available.

In the special case where all the action sets are singletons (and thus the Controller is a dummy), Kolmogorov’s 0-1 law says that the probability for the goal to be reached is 0 or 1. We construct a number of counterexamples to show that in general the value of the control problem can be strictly between 0 and 1, and derive several sufficient conditions for the 0-1 ``law" to hold.
\end{abstract}

\noindent\textbf{Keywords:} Tail set, zero-one law, Markov decision process, maximal branching process, branching process in varying environments.\\

\noindent\textbf{MSC:} 
60F20; 
90C40; 
60J85.

\section{Introduction}
\noindent\textbf{The framework:} The Controller (C) chooses an infinite sequence of actions, $a_{0},a_{1},\ldots$, one at a time. Having taken actions $h_{t} = (a_{0},\ldots,a_{t-1})$ prior to stage $t$, the Controller's action $a_{t}$ at stage $t$ is required to be an element of the non-empty, finite subset $A(h_{t})$ of $\N$. The Controller's goal is to maximize the probability that the sequence $(a_0,a_1,\ldots)$ is an element of a given subset $G$ of $\N^\N$. The set $G$, called the \textit{goal}, is assumed to be a Borel tail set. The set $A(h_{t})$ is chosen from a distribution $p_{t}$, called the \textit{primitive distribution}, independently over all $t \in \N$ and all $h_{t} \in \N^{t}$.

One can equivalently represent the restriction on the Controller's choices with the help of $\omega$, a finitely branching pruned tree on $\N$, called the \textit{decision tree}. The Controller is required to choose the actions so that each finite sequence $h_{t}$ of actions is an element of $\omega$. The action set $A_{\omega}(h_{t})$ is then defined as consisting of actions $a_{t} \in \N$ such that $h_{t}a_{t} \in \omega$. This is the formalism we choose to follow in the main body of the paper. We define a measure on the set of all decision trees in the spirit of the familiar Galton-Watson measure (Neveu \cite{Neveu86}, Lyons and Peres \cite{LyonsPeres17}). Its key feature is that, given the event $\{h_{t} \in \omega\}$, the subtrees following the sequence $h_{t}$ are conditionally independent, and the action set $A_{\omega}(h_{t})$ is distributed according to $p_{t}$.

We study two information structures, defined by how much of the decision tree is revealed to the Controller before she takes each action.\medskip

\noindent\textbf{Extending Kolmogorov's zero-one law:} A special case of our model is one where the primitive distributions $p_{t}$ only assign a positive probability to singleton action sets. In this case, the Controller is a dummy, and how much of the decision tree she gets to know is irrelevant. Kolmogorov's zero-one law says that in this special case, the probability for the goal to be reached is either $0$ or $1$. The aim of this paper is to extend this result beyond the special case.\medskip

\noindent\textbf{The omniscient and the $m$-foresight information structures:} Under the \textit{omniscient} information structure, serving as a benchmark, the decision tree is entirely revealed to the Controller before the very first action she has to take. The omniscient Controller is thus able to condition her actions at any stage on the entire decision tree.  

In contrast, under $m$-\textit{foresight}, the decision tree is only gradually being revealed to the Controller in the process of decision-making, $m$ being the depth of foresight. 

Thus, under $0$-foresight, at a stage $t$ the Controller is only revealed the set of actions available to her at that stage. Once the action at stage $t$ is chosen, the stage $t+1$ set of actions is revealed to her, and so on. Under $1$-foresight, the Controller not only knows what actions are available to her at the current stage, but also what actions will or will not be available at the next stage. More generally, under the scenario of $m$-foresight, at any given stage, the next $m$ generations of the decision tree are revealed to the Controller before she chooses an action.\medskip

\noindent\textbf{The omniscient and the $m$-foresight values:} The omniscient value is defined as the probability of the event that the decision tree $\omega$ has an infinite branch that is an element of $G$. The definition is motivated as follows: if the body of the decision tree (i.e. the set of its infinite branches) is disjoint from the goal, then the goal cannot be reached, not even by the omniscient Controller. If, however, they do have a point in common, then the omniscient Controller, who knows the entire decision tree, can simply execute any sequence of actions that is an element of the goal. One can also give a more decision-theoretic definition of the omniscient value as the supremum of the probability of success over all omniscient strategies. 

The Controller with $m$-foresight does not have complete knowledge of the decision tree, and she may only condition her decision at stage $t$ on the information revealed to her up to that stage. This places a restriction on the set of strategies available to the Controller. The $m$-foresight value is defined as the supremum of the probability of success over all $m$-foresight strategies. 

The $m$-foresight value of a control problem (CP) is non-decreasing in $m$, and the omniscient value is at least as large as the $m$-foresight value for each $m$.\medskip

\noindent\textbf{Results:} The $0$-foresight value of any control problem is either $0$ or $1$ (Corollary \ref{thm.0foresight}). In contrast, neither the omniscient value nor the $m$-foresight value, for $m \geq 1$, obey the zero-one “law”. We construct a number of counterexamples to this effect, and supply several sufficient conditions for the law to hold. 

The first sufficient condition requires a certain time-invariance of the control problem. More precisely, if the primitive distributions are time-invariant and the goal is  shift-invariant, then both the omniscient and the $m$-foresight values, for each $m \in \N$, are $0$ or $1$ (Theorems \ref{thm.omni:shift} and \ref{thm.finite:shift}). We note that the condition of shift-invariance can in fact be relaxed. Neither the assumption of shift-invariance of the goal, nor the assumption of  time-invariance of the primitive distributions, on their own, would be sufficient (Examples \ref{exl.1foresight:shift} and \ref{exl.1foresight:invariant}).

The second sufficient condition, in contrast to the first one, does not place any restrictions on the goal, and allows for a time-dependence of primitive distributions. It requires instead that the cardinality of the decision set at all stages $t$ be dominated by a single random variable with finite mean. Formally, let $\#p_{t}$ be the distribution of the cardinality of the set of actions at stage $t$. If the distributions $\#p_{t}$ are all dominated by a single distribution with finite mean, then the omniscient (resp. any $m$-foresight) values are equal to $0$ or $1$ (Theorems \ref{thm.omni:dominance} and \ref{thm.all:dominance}). Example \ref{exl.lyons} illustrates the tightness of the result.

The conditions discussed above are sufficiently strong to ensure that the omniscient and all the $m$-foresight values obey the zero-one law. Focusing on the $m$-foresight value for a given $m \in \N$, one can obtain a weaker sufficient condition. This is accomplished by Theorem \ref{thm.finite:dominance}. The condition of Theorem \ref{thm.finite:dominance} is linked to the notion of maximal branching process (Lamperti \cite{Lamperti70,Lamperti72}). Example \ref{exl.lyons} illustrates that the condition of Theorem \ref{thm.finite:dominance} does not guarantee the zero-one law for the omniscient value.\medskip

\noindent\textbf{Techniques:} The proofs of the results for the omniscient value rely on techniques familiar from the literature on branching process (e.g. Athreya and Ney \cite{AthreyaNey04}, Lyons and Peres \cite[Chapter 5]{LyonsPeres17}). The proofs of the results pertaining to the $m$-foresight have a rather different nature. We cast this scenario as a setting of a Markov decision process (MDP) and use some known facts about the value of an MDP (e.g. Bertsekas and Shreve \cite{BertsekasShreve96}). The proof also employs coupling of two stochastic processes (Lindvall \cite{Lindvall02}). As noted earlier, the sufficient condition for $m$-foresight makes use of the notion of a maximal branching process.\medskip        
 
\noindent\textbf{Related literature:} The measure on the space of decision trees we employ is similar to the classical Galton-Watson measure on the space of Galton-Watson trees (Neveu \cite{Neveu86}, Lyons and Peres \cite{LyonsPeres17}). Moreover, the size of the $t$-th generation of a decision tree is distributed like the branching process in varying environments (e.g. Fearn \cite{Fearn71}).  

The scenario of finite foresight is intimately linked to maximal branching processes (Lamperti \cite{Lamperti70, Lamperti72}). In the situation of 1-foresight, the distribution of the number of actions that the Controller has at her disposal is dominated by a certain maximal branching process. It is this property that drives Theorem \ref{thm.finite:dominance}. 

Decision problems and games where the Controller's goal is a tail set are ubiquitous in computer science. The classical B\"{u}chi, co-B\"{u}chi, parity, Muller, Streett, and Rabin objectives all have the tail property (we refer to Apt and Gr\"{a}del \cite{AptGradel11} and Chatterjee and Henzinger \cite{Chatterjee12} for an overview). The existence of equilibria in multiplayer games with tail objective functions has been studied in Ashkenazi-Golan, Flesch, Predtetchinski, and Solan \cite{AFPS22}, and Flesch and Solan \cite{FleschSolan23a,FleschSolan23b}. 

There are related works on random dynamic games (Basu, Holroyd, Martin, and Wästlund \cite{Basu16}, Flesch, Predtetchinski, and Suomala \cite{FleschSuomala23}, Garnier and Ziliotto \cite{Garnier22}, Holroyd and Martin \cite{HolroydMartin21}, Holroyd, Marcovici, and Martin \cite{HolroydMarcoviciMartin19}): in these models two antagonistic Controllers interact in an environment randomly chosen by Nature; at the focus of these studies is the distribution of the value of the random game. 

Markov decision processes were introduced by Shapley \cite{Shapley53} and Blackwell \cite{Blackwell67}. In contrast to much of the literature on the subject that considers payoff functions defined by aggregating stage rewards (such as discounted, total, average reward evaluations), Markov decision process considered here has an objective represented by a Borel tail set. General payoff functions in Markov decision processes  have been considered by Feinberg \cite{Feinberg83}. Gambling house (Dubins, Savage, Sudderth, Gilat \cite{DubinsSavageSudderthGilat14} and Dubins, Maitra, and Sudderth \cite{DubinsMaitraSudderth02}) is another example of a Markov decision process with a Borel-measurable (or upper semi-analytic) objective function. 

Economists have considered the idea that the decision maker can only see a finite number of stages of the decision tree, see for example Ke \cite{Ke19}, and the recent paper by Pycia and Troyan \cite{PyciaTroyan23}, and the references therein. In particular, Pycia and Troyan focus on the scenarios of omniscience, 1-foresight, and 0-foresight. Unlike much of the economic literature on the subject that treats bounded foresight as a form of bounded rationality, our approach is that of bounded information, and is thus perfectly compatible with the standard game/decision theoretic paradigm of rational decision-making. 

The paper is structured as follows. In Section 2 we introduce control problems with random action sets. In Section 3 we define the omniscient value and the $m$--foresight value. Section 4 presents the main results of the paper and a number of examples. Section 5 to Section 7 are dedicated to the proofs. Section 8 concludes with some additional remarks.

\begin{table}\vspace{-2cm}
\renewcommand{\arraystretch}{1.2}
\begin{center}
\begin{tabular}[h]{lll}\hline\toprule
notation & description & page\\\midrule
$\A$ & collection of non-empty finite subsets of $\N$ & \pageref{nota.A}\\
$p_{t}$ & primitive distribution at stage $t$ & \pageref{nota.pt}\\
$G$ & objective & \pageref{nota.G}\\
$\omega$ & decision tree & \pageref{nota.omega}\\
$A_{\omega}(h)$ & action set at $h$ & \pageref{nota.A(h)}\\
$\omega(h)$ & subtree of $\omega$ at $h$ & \pageref{nota.omega(h)}\\
$\omega_{n}$ & generation $n$ of $\omega$ & \pageref{nota.omega_n}\\
$[\omega]$ & body of $\omega$ & \pageref{nota.[omega]}\\
$\Omega$ & set of decision trees & \pageref{nota.Omega}\\
$\F_{t}$ & sigma-algebra generated by $\omega_{t}$ & \pageref{nota.Ft}\\
$\mu$, $\mu_{t}$ & measures on $\Omega$ & \pageref{nota.mu}\\
$\#p_{t}$ & distribution of cardinality of stage $t$ action set & \pageref{nota.cardpt}\\
$\#\omega_{m}$ & cardinality of $\omega_{m}$ & \pageref{nota.cardomega_n}\\
$\#p_{t}^{m}$ & distribution of cardinality of $\#\omega_{m}$ under $\mu_{t}$& \pageref{nota.cardptm}\\
$\Omega_{G}$ & event of success of the omniscient Controller & \pageref{nota.Omega_G}\\
$v^{\infty}$ & omniscient value & \pageref{nota.vinfty}\\
$\sigma$ & strategy & \pageref{nota.sigma}\\
$\Omega_{G}^{\sigma}$ & event of success of strategy $\sigma$ & \pageref{nota.Omega_G^sigma}\\
$v^{m}$ & $m$-foresight value & \pageref{nota.v^m}\\
$G_{t}$ & $t$-shifted goal & \pageref{nota.Gt}\\
$\phi_{t}$ & unconditional distribution of state at stage $t$ & \pageref{nota.phit}\\
$Z_{t}^{\epsilon}$, $z_{t}^{\epsilon}$ & set and probability of good states & \pageref{nota.zt}\\
MDP$^{m}$ & MDP for $m$-foresight & \pageref{nota.MDP}\\
$S = S_{m}$ & state space of MDP$^{m}$ & \pageref{nota.S}\\
$A(s)$ & action set of MDP$^{m}$ & \pageref{nota.A(s)}\\
$\Omega^{*}$ & measure space of MDP$^{m}$ & \pageref{nota.Omega*}\\
$\I_{t}$ & stage $t$ sigma-algebra on $\Omega^{*}$ & \pageref{nota.It}\\
$G^{*}$ & the goal in MDP$^{m}$ & \pageref{nota.G*}\\
$\mu_{\sigma}^{*}$ & measure on $\Omega^{*}$ induced by $\sigma$ & \pageref{nota.mu*}\\
$v_{t}^{m}$ & stage $t$ value of MDP$^{m}$ & \pageref{nota.vtm}\\
$s_{t}^{\sigma}$ & state at stage $t$ encountered under $\sigma$ & \pageref{nota.stsigma}\\
$u(s)$ & size of state $s$ & \pageref{nota.size}\\\bottomrule
\end{tabular}
\caption{Table of notation}
\end{center}
\end{table}

\section{Control problems with random action sets}
\noindent\textbf{Control problems with random action sets:} Let $\A$ \label{nota.A} be the collection of all non-empty finite subsets of $\N = \{0,1,\ldots\}$. A \textit{control problem with random action sets} (CP) is a tuple $(p_{0},p_{1},\ldots,G)$ where for each $t \in \N$, $p_{t}$ \label{nota.pt} is a probability distribution on $\A$, and $G$ \label{nota.G} is a subset of $\N^\N$. The probability distribution $p_{t}$ is called the \textit{stage $t$ primitive distribution}. The set $G$ is called the \textit{goal}. Throughout the paper $G$ is assumed to be a Borel tail set. The latter means that for each $x$ and $y$ in $\N^\N$ such that $x_{n} = y_{n}$ for all but finitely many $n \in \N$, one has $x \in G$ if and only if $y \in G$.\medskip

\noindent\textbf{Decision trees:} Let $\N^* = \cup_{n \in \N}\N^n$ be the set of all finite sequences of natural numbers, including the empty sequence $\o$. By a \textit{decision tree} we mean a pruned finitely splitting tree on $\N$. More explicitly, a decision tree is a subset $\omega$ \label{nota.omega} of $\N^*$ satisfying the following conditions: (a) it is closed under prefixes: whenever $h \in \N^*$ and $a \in A$ are such that the sequence $ha$ is in $\omega$, so is the sequence $h$; (b) for each $h \in \omega$ there exists at least one but at most finitely many $a \in \N$ such that $ha \in \omega$. Elements of a decision tree are called its nodes. For a node $h \in \omega$ the \textit{action set} at $h$ is $A_{\omega}(h) = \{a \in \N:ha \in \omega\}$\label{nota.A(h)}. 

Given $h \in \omega$, the subtree of $\omega$ at $h$ is the decision tree $\omega(h)$ \label{nota.omega(h)} consisting of the sequences $g \in \N^*$ such that $hg \in \omega$. For $n \in \N$ let $\omega_{n}$ \label{nota.omega_n} denote the set $\omega \bigcap \N^n$; this is the set of nodes of $\omega$ of length exactly $n$, and is called the $n$th \textit{generation} of the tree; in particular, $\omega_{0} = \{\o\}$. Note that, since the decision tree is closed under prefixes, the $n$th generation is completely determined by its $(n+1)$st generation. Let $[\omega]$ \label{nota.[omega]} be the body of $\omega$, that is, the set of $x \in \N^\N$ such that $(x_{0},\ldots,x_{n}) \in \omega$ for each $n \in \N$.\medskip

\noindent\textbf{The topology on the set of decision trees:} Let $\Omega$ \label{nota.Omega} be the set of all decision trees. We introduce the topology on $\Omega$ with the help of the metric $d$ defined as follows: for two distinct decision trees $\omega$ and $\eta$ let $d(\omega,\eta) = 2^{-t}$ where $t \in \N$ is the smallest natural number such that $\omega_{t} \neq \eta_{t}$. This makes $\Omega$ a Polish space. In fact, the bijection $\omega \mapsto [\omega]$ from $\Omega$ to $\mathcal{K}(\N^\N)$ (Kechris \cite[Proposition 2.4 and Excercise 4.11]{Kechris95}), the space of non-empty compact subsets of $\N^\N$ with the Vietoris topology, is a homeomorphism.

Let $\mathscr{F}_{t}$ \label{nota.Ft} denote the sigma-algebra on $\Omega$ generated by the function $\Omega \to \N^t$, $\omega \mapsto \omega_{t}$ (in particular, $\mathscr{F}_{0}$ is the trivial sigma-algebra $\{\o,\Omega\}$). Let $\mathscr{F}$ be the sigma-algebra generated by the set $\cup_{t \in \N}\mathscr{F}_{t}$. This is exactly the Borel sigma-algebra on $\Omega$.\medskip

\noindent\textbf{The measure on $\Omega$:} Given a CP as above, we define a Borel probability measure $\mu = p_{0} \sqcap p_{1} \sqcap \cdots$ \label{nota.mu} on $\Omega$ as follows:  given an $\omega \in \Omega$ and $t \in \N$ let
\[\mu(\{\eta \in \Omega: \eta_{t} = \omega_{t}\}) = \prod_{k < t}\prod_{h \in \omega_{k}} p_{k}(A_{\omega}(h)).\]
This specifies $\mu$ on the sigma-algebra $\mathscr{F}_{t}$ for each $t \in \N$, and hence, by Kolmogorov's extension theorem, on $\mathscr{F}$. We also define $\mu_{t} = p_{t} \sqcap p_{t+1} \sqcap \cdots$. \label{nota.mu_t} In particular, $\mu_{0} = \mu$.

These measures have the following properties:
\begin{enumerate}
\item Under $\mu$, the random variable $\omega \mapsto \omega_{1}$ is distributed according to $p_{0}$, that is: $\mu(\{\omega_{1} = A\}) = p_{0}(A)$ for each non-empty finite subset $A$ of $\N$. 
\item For each stage $t \in \N$, under the conditional measure $\mu(\cdot|\mathscr{F}_{t})$ the random variables $(\omega(h):h \in \omega_{t})$ are independent and each is distributed according to $\mu_{t}$.  
\end{enumerate}

Unless stated otherwise, $\mathbb{E}$ denotes the expectation with respect to the measure $\mu_{0}$. With a slight abuse of notation we write $\mu_{0}$ to denote the completion of $\mu_{0}$.\medskip

\noindent\textbf{The process $\#\omega_{t}$ and related notions:} Let us denote by $\#p_{t}$ \label{nota.cardpt} the distribution of the cardinality of the action set under the measure $p_{t}$: thus $\#p_{t}(n) = p_{t}(\{A \subset \N:\#A = n\})$ for each $n \in \{1,2,\ldots\}$. Let $g_{t}$ \label{nota.gt} be the probability generating function of $\#p_{t}$: $g_{t}(x) = \#p_{t}(1) \cdot x + \#p_{t}(2) \cdot x^{2} + \cdots$. For $m \in \N$ define $g_{t}^{m+1} = g_{t} \circ \cdots \circ g_{t+m}$ \label{nota.gtm} (and let $g_{t}^{0}$ be the identity map); let $\#p_{t}^{m}$ \label{nota.cardptm} denote the distribution whose probability generating function is $g_{t}^{m}$ (so in particular, $\#p_{t}^{0}$ is the Dirac measure on 1, and $\#p_{t}^{1} = \#p_{t}$). 

Denote by $\#\omega_{n}$ \label{nota.cardomega_n} the cardinality of the $n$th generation of $\omega$. Fix a $t \in \N$. One can see that the process $\#\omega_{0}, \#\omega_{1}, \ldots$ on $(\Omega, \mathscr{F}, \mu_{t})$ is the branching process with generation dependence (also known as branching process in a varying environment, BPVE). We have $\#p_{t}^{m}(n) = \mu_{t}(\{\omega \in \Omega: \#\omega_{m} = n\})$ for each $n \in \{1,2,\ldots\}$ (Fearn \cite[Proposition 1]{Fearn71}). That is, $\#p_{t}^{m}$ is the distribution of $\#\omega_{n}$ under the measure $\mu_{t}$.\medskip

\section{The omniscient and the $m$-foresight values}
We investigate two information structures differing in the information that the Controller has at her disposal, and introduce the corresponding values.

\subsection{The omniscient value}
Let $\Omega_{G} = \{\omega \in \Omega: [\omega] \cap G \neq \o\}$\label{nota.Omega_G}. The event $\Omega_{G}$ is an analytic (see the Appendix for details) subset of $\Omega$. The \textit{omniscient value} of the CP $(p_{0},p_{1},\ldots,G)$ is defined as \label{nota.vinfty}
\[v^{\infty}(p_{0},p_{1},\ldots,G) = \mu_{0}(\Omega_{G}).\]

The motivation for the definition is as follows. Suppose that Nature first selects the decision tree $\omega \in \Omega$ from the measure $\mu_{0}$, and then reveals it to the Controller. Having observed the entire decision tree $\omega$, the Controller proceeds with making the decisions. We refer to this as the scenario of an omniscient Controller. 

Clearly then, the omniscient Controller is able to reach her goal as soon as the body $[\omega]$ of the decision tree $\omega$ has a point in common with the objective $G$. On the other hand, if the two sets are disjoint, the goal cannot be reached at all. Thus $\Omega_{G}$ could be interpreted as the event that the omniscient Controller succeeds in reaching the goal, and the omniscient value as the probability of success of the omniscient Controller.  

One can also give an equivalent definition of the omniscient value more in line with decision-theoretic tradition. Even though it is not needed to derive our results, we state this equivalent definition in anticipation of the discussion of $m$-foresight in the following subsection.      

Let $\sigma = (\sigma_{0},\sigma_{1},\ldots): \Omega \to \N^\N$\label{nota.sigma} be an analytically measurable function (i.e. a function measurable with respect to the sigma-algebra generated by analytic sets). It is called a \textit{strategy} if $\sigma(\omega) \in [\omega]$ for each $\omega \in \Omega$. Equivalently, $\sigma$ is a strategy if $(\sigma_{0}(\omega), \ldots, \sigma_{n}(\omega)) \in \omega$ for each $\omega \in \Omega$ and $n \in \N$. Given a strategy $\sigma$ let us also define $h_{n}^{\sigma}(\omega) = (\sigma_{0}(\omega), \ldots, \sigma_{n-1}(\omega))$, with $h_{0}^{\sigma}(\omega) = \o$, to be the node of the decision tree $\omega \in \Omega$ reached by the Controller at stage $n$, provided she adheres to the strategy $\sigma$.

Let $\Sigma_{\infty}$ be the set of strategies. For $\sigma \in \Sigma_{\infty}$ we define $\Omega_{G}^{\sigma} = \{\omega \in \Omega: \sigma(\omega) \in G\}$\label{nota.Omega_G^sigma}. This is the event that the strategy $\sigma$ is successful in reaching the goal $G$. The omniscient value satisfies
\[v^{\infty}(p_{0},p_{1},\ldots,G) = \sup_{\sigma \in \Sigma_{\infty}} \mu_{0}(\Omega_{G}^{\sigma}),\]
where the supremum is attained (see Appendix for details). This equation could be taken as a game-theoretic definition of the omniscient value.  

\subsection{The $m$-foresight value}
This scenario proceeds as follows. Nature first draws the decision tree $\omega$ from the measure $\mu_{0}$, and reveals it to the Controller gradually, as she proceeds with making her decisions. At any stage of the decision-making process, the Controller only knows what decision nodes she will be able to reach in $m+1$ stages. We can thus represent the Controller's problem with $m$-foresight as follows\footnote{With a slight abuse of notation we write $\omega_{m+1}(h)$, rather than $(\omega(h))_{m+1}$, to denote the $(m+1)$st generation of the subtree $\omega(h)$.}:
\begin{center}
\begin{tabular}{rccccccc}
{\rm Nature} & $\omega_{m+1}(\o)$ && $\omega_{m+1}(a_{0})$ && $\omega_{m+1}(a_{0},a_{1})$ && $\cdots$\\ {\rm Controller} && $a_{0}$ && $a_{1}$ && $a_{2}$ & $\cdots$
\end{tabular}
\end{center}

Thus, prior to the Controller making her first choice, Nature reveals the first $(m+1)$ generations of the decision tree. Once the Controller makes her first decision, $a_{0}$, she is revealed the $(m+1)$st generation of the subtree $\omega(a_{0})$. The Controller then takes her second action, $a_{1}$, and she is revealed the $(m+1)$st generation of the subtree $\omega(a_{0},a_{1})$, and so on. The Controller's action should satisfy the condition $a_{t} \in A_{\omega}(a_{0},\ldots,a_{t-1})$ for each $t \in \N$. 

We give a formal definition of $m$-foresight strategies. A strategy $\sigma$ has $m$-\textit{foresight} if it meets the following condition: Let $\omega$ and $\eta$ be decision trees, let $x = \sigma(\omega)$ and $y = \sigma(\eta)$, and $t \in \N$. If $\omega_{m+1}(\o) = \eta_{m+1}(\o)$, then $x_{0} = y_{0}$; if in addition $\omega_{m+1}(x_{0}) = \eta_{m+1}(y_{0})$, then $x_{1} = y_{1}$; if in addition $\omega_{m+1}(x_{0},x_{1}) = \eta_{m+1}(y_{0},y_{1})$, then $x_{2} = y_{2}$; and so on. The definition thus requires that the actions chosen by the strategy at a particular stage only depend on the information revealed to the Controller up to that stage. Let $\Sigma_{m}$ be the set of strategies with $m$-foresight.

\begin{exl}\label{exl.benchmark}\rm
An example of a $0$-foresight strategy is one that always picks the smallest action in the action set. It is (recursively) defined by letting $\sigma_{t}(\omega)$ be the smallest element in the set $A_{\omega}(h_{t}^{\sigma}(\omega))$.

Suppose instead that the Controller chooses the action at stage $t$ so as to maximize the number of next-stage actions. Formally, $\sigma_{t}(\omega)$ is the smallest action $a_{t} \in A_{\omega}(h_{t}^{\sigma}(\omega))$ maximizing $\#A_{\omega}(h_{t}^{\sigma}(\omega),a_{t})$. This is an example of a 1-foresight strategy, which we refer to as \textit{1-step maximizing}. $\Box$
\end{exl}

The $m$-foresight value of the CP, denoted by $v^m$, is defined as the supremum of the success probabilities under $m$-foresight strategies:
\[v^{m}(p_{0},p_{1},\ldots,G) = \sup_{\sigma \in \Sigma_{m}} \mu_{0}(\Omega_{G}^{\sigma}).\]\label{nota.v^m}

It is clear that for any CP, the $m$-foresight value is non-decreasing in $m \in \N$, and is not larger than the omniscient value: $v^{0}  \leq v^{1} \leq v^{2} \leq \cdots \leq v^{\infty}$.

\section{Results}
We start out by stating two facts that apply to all CPs. The first says that the $0$-foresight value obeys the zero-one law. 

\begin{prop}\label{thm.0foresight}
It holds that $v^{0}(p_{0},p_{1},\ldots,G) \in \{0,1\}$. 
\end{prop}

\noindent This fact has a curious implication that likewise applies to all CPs: either the (omniscient/ $m$-farsighted) C succeeds almost surely for the goal $G$, or she succeeds almost surely for the goal $\N^{\N} \setminus G = G^{c}$.

\begin{cor}\label{thm.01lawI} 
Let $m \in \N \cup \{\infty\}$. Then $v^{m}(p_{0},p_{1},\ldots,G) = 1$ or $v^{m}(p_{0},p_{1},\ldots,G^c) = 1$. 
\end{cor}
\begin{proof}
If $v^{0}(p_{0},p_{1},\ldots,G) = 0$ then $v^{m}(p_{0},p_{1},\ldots,G) = 1$. Suppose that $v^{0}(p_{0},p_{1},\ldots,G) = 1$ and take any 0-foresight strategy $\sigma$. Then $\mu_{0}(\Omega_{G}^{\sigma}) = 0$ and hence $\mu_{0}(\Omega_{G^c}^{\sigma}) = 1$. Hence $v^{0}(p_{0},p_{1},\ldots,G^c) = 1$ and $v^{m}(p_{0},p_{1},\ldots,G^c) = 1$.
\end{proof}

Unlike the $0$-foresight value, neither the $m$-foresight, for $m > 0$, nor the omniscient value obeys the zero-one law. We shall see three counterexamples to this effect. Presently we discuss one where both the $m$-foresight, for $m > 0$, and the omniscient values are neither $0$ nor $1$. 

\begin{exl}\rm\label{exl.1foresight:shift}
For $t \in \N$ let $p_t$ be the measure on $\A$ that assigns probability to only two sets: the probability $1-2^{-t}$ to the set $\{0\}$ and the probability $2^{-t}$ to the set $\{0,1,\ldots,t2^{t}\}$ (Table \ref{table.1foresight:shift}).

\begin{table}
\begin{center}
\renewcommand{\arraystretch}{1.4}
\begin{tabular}{lc}\toprule
set & probability\\\midrule
$\{0\}$ & $1 - 2^{-t}$\\
$\{0,1,\ldots,t2^{t}\}$ & $2^{-t}$\\\bottomrule
\end{tabular}
\end{center}
\caption{The primitive distribution $p_{t}$ in Example \ref{exl.1foresight:shift}}\label{table.1foresight:shift}
\end{table}

Let $G = \{x \in \N^\N: x_t = 0\text{ for at most finitely many }t \in \N\}$. We argue that 
\[0 < v^{1}(p_{0},p_{1},\ldots,G) \leq v^m(p_{0},p_{1},\ldots,G) \leq v^\infty(p_{0},p_{1},\ldots,G) < 1.\]
Clearly, it suffices to prove the first and the last inequality.\medskip

\noindent\textsc{Claim 1:} $v^\infty(p_{0},p_{1},\ldots,G) < 1$, Let $\delta$ be the decision tree consisting of the sequences of $0$'s. Of course, $[\delta]$ is a singleton set containing the path $000\cdots$ only, which is not in $G$. Thus the tree $\delta$ is not in $\Omega_{G}$, and consequently it suffices to argue that $\mu(\{\delta\}) > 0$. And indeed, $\mu(\{\delta\}) = \prod_{t \in \N}(1-2^{-t-1}) > 0$.\medskip

\noindent\textsc{Claim 2:} $0 < v^{1}(p_{0},p_{1},\ldots,G)$. Let $\sigma$ be the strategy that, at any stage $t$, chooses the smallest non-zero action that could be followed up by a non-zero action at stage $t+1$; if all actions at stage $t$ have a single follow-up action (action $0$) at stage $t+1$, pick action $1$. We argue that, under the strategy $\sigma$, there is a positive probability that the Controller never takes action $0$, except at stage $0$.

Let $E_{0} = \Omega$ and for $t \geq 1$ let $E_{t}$ be the event that at stage $t$ the Controller has a non-zero action at her disposal. Note that, in the event $E_{t}$, the Controller does take a non-zero action at stage $t$. We show that the event $\cap_{t \geq 1}E_{t}$ has positive probability. Clearly, $\mu(E_{1} | E_{0}) = \mu(E_{1}) = \tfrac{1}{2}$. For $t \geq 1$, conditional on the event $E_{t} \cap \cdots \cap E_{0}$, the event $E_{t+1}$ does not occur exactly when each of the $t2^{t}$ non-zero actions at stage $t$ have only action $0$ as a possible follow up. This happens with probability $(1-2^{-t-1})^{t2^{t}}$. Since $(1-2^{-t-1})^{2^{t}}$ approaches $e^{-1/2}$ as $t \to \infty$, we have
\[\sum_{t \geq 1}(1-2^{-t-1})^{t2^{t}} < \infty.\]
Hence
\[\mu(\cap_{t \geq 1}E_{t}) = \prod_{t \geq 1}\mu(E_{t+1} | E_{t} \cap \cdots \cap E_{0}) = \prod_{t \geq 1}[1-(1-2^{-t-1})^{t2^{t}}] > 0.\,\Box\]
\end{exl}

We turn to sufficient conditions for the zero-one law to hold. The first sufficient condition involves a certain type of time-invariance of the CP. Let us say that the primitive distribution is \emph{time-invariant} if $p = p_{0} = p_{1} = \cdots$. For every $t \in \N$, we define the $t$-shifted goal $G_{t}$\label{nota.Gt} as the set of sequences $x \in \N^\N$ such that $(0,\ldots,0,x_{0},x_{1},\ldots) \in G$ (with $t$ zeros). We remark that, since $G$ is a tail set, the sequence $x \in \N^\N$ is an element of the $t$-shifted goal $G_{t}$ if and only if $(a_{0}, \ldots, a_{t-1}, x_{0}, x_{1}, \ldots,) \in G$ for some sequence $(a_{0}, \ldots, a_{t-1}) \in \N^{t}$ of length $t$. In particular, the $0$-shifted goal $G_0$ is equal to the original goal $G$. 

For the ease of referencing, we present the results on the omniscient value and the $m$-foresight value as two separate theorems, even if the sufficient condition is the same in both.  

\begin{thm}\label{thm.omni:shift}
Suppose that the primitive distribution is time-invariant and that there exist $0 \leq k_0 < k_1$ such that $G_{k_{0}} \subseteq G_{k_{1}}$. Then $v^\infty(p,p,\ldots,G) \in \{0,1\}$.
\end{thm}

\begin{thm}\label{thm.finite:shift}
Suppose that the primitive distribution is time-invariant and that there exist $0 \leq k_0 < k_1$ such that $G_{k_{0}} \subseteq G_{k_{1}}$. Then $v^{m}(p,p,\ldots,G) \in \{0,1\}$ for each $m \in \N$. 
\end{thm}

The condition the two theorems impose on the goal could be interpreted as follows: starting from the stage $k_{0}$, the goal becomes more permissive with the passing of every $k_{1} - k_{0}$ stages. And thus, beginning from the stage $k_{0}$, the Controller will be periodically facing the same or even a “simpler” CP. The condition is in particular satisfied if the goal is shift-invariant, i.e. if $G = G_{1}$ (in which case $G = G_{t}$ for every $t \in \N$). Let us highlight that the two theorems do not impose any restriction on the primitive distribution apart from time-invariance; in particular, the mean number of actions may be infinite.

Example \ref{exl.1foresight:shift} above shows that shift-invariance of the goal, without the condition of time-invariance of the primitive distribution, is not sufficient for obtaining the zero-one law. Example  \ref{exl.1foresight:invariant} below illustrates that likewise time-invariance of the primitive distribution, by itself, is not sufficient either. Note that in this example, the $t$-shifted goals form a strictly decreasing sequence, and thus the condition of Theorems \ref{thm.omni:shift} and \ref{thm.finite:shift} is not satisfied. 

\begin{exl}\label{exl.1foresight:invariant}\rm
Let $m_{t} \in \N$ and $r_{t} \in (0,1)$ be any sequences such that 
\begin{itemize}
\item[(a)] $1 = m_{0} < m_{1} < \cdots$ and $0 < r_{0} < r_{1} < \cdots$ with $r_{t}$ approaching $1$,
\item[(b)] $\prod_{t \geq 1}r_{t}^{m_{0} \cdot \cdots \cdot m_{t-1}} > 0$, and
\item[(c)] $\prod_{t \in \N}(1 - r_{t}^{m_{t}}) > 0$ 
\end{itemize}
For example, one can define the sequence $\{m_{t}\}_{t \in \N}$ recursively by $m_{0} = 1$ and $m_{t} =  t^{3} \cdot \prod_{i < t}m_{i}$ for each $t \geq 1$, and let $r_{t} = 2^{-t/m_{t}}$ for $t \in \N$.

Let $\{M_{t}:t \in \N\}$ be a partition of $\N$ such that the cardinality of $M_{t}$ is $m_{t}$. 
Let $p = p_{0} = p_{1} = \cdots$ be the time-invariant primitive distribution that only assigns positive probability to the action sets in the collection $\{M_{t}:t \in \N\}$, and such that the probability for $p$ to choose an action set from the subcollection $\{M_{0}, \ldots, M_{t}\}$ equals $r_{t}$. Let the goal be given by 
\[G = \bigcup_{T \in \N}\bigcap_{t \geq T} \{x \in \N^\N: x_{t} \in M_{t+1}\cup M_{t+2} \cup \cdots\}.\]
We claim that 
\[0 < v^{1}(p_{0},p_{1},\ldots,G) \leq v^m(p_{0},p_{1},\ldots,G) \leq v^\infty(p_{0},p_{1},\ldots,G) < 1.\]

We relegate the rigorous proof of the claims to the next section, and give here only a gist of the arguments.

To show that the omniscient value is smaller than 1, we argue that, with positive probability, at no stage $t \in \N$ is there a node of the decision tree with more than $m_{t}$ actions. This means that there is a positive probability for the decision tree $\omega$ to be such that $x_{t} \in M_{0} \cup \cdots \cup M_{t}$ for each infinite branch $x$ of $\omega$ and each stage $t \in \N$. Indeed, with probability $r_{0}$ there are only $m_{0}$ actions at the root of the tree. Conditional on this event, with probability $r_{1}^{m_{0}}$ there are at most $m_{1}$ actions at each of the nodes of generation 1. Conditional on these two events, generation 2 of the tree has no more than $m_{0}m_{1}$ nodes, and hence (still conditional on these events), with probability $r_{2}^{m_{0}m_{1}}$ there are at most $m_{2}$ actions at each of the nodes of generation 2. Continuing in this way, we find that, with probability $r_{0} \cdot r_{1}^{m_{0}} \cdot r_{2}^{m_{0}m_{1}} \cdots r_{t}^{m_{0} \cdots m_{t-1}} \cdots > 0$, at no stage $t \in \N$ is there a node with more than $m_{t}$ actions.   

To show that the $1$-foresight value is positive, consider a strategy $\sigma$ that maximizes the number of follow-up, or next-stage, actions available; in other words, a strategy $\sigma$ that picks, at any stage $t \in \N$, an action that leads to a node with the largest action set. We show that under $\sigma$, with positive probability, at any given stage $t$ the Controller has at least $m_{t+1}$ actions available. 

Indeed, with probability $1-r_{0}$, there are at least $m_{1}$ actions at stage $0$. Conditional on this event, with probability $1 - r_{1}^{m_{1}}$, at least one of the nodes in generation 1 will have no fewer than $m_{2}$ actions. In this event, in view of the definition of $\sigma$, the Controller will have at least $m_{2}$ actions at stage $1$. Conditional on the latter event, with probability $1 - r_{2}^{m_{2}}$ the Controller will have at least $m_{3}$ actions at stage $2$, and so on. And thus we find that, with probability  $(1-r_{0}) \cdot (1 - r_{1}^{m_{1}}) \cdot (1 - r_{2}^{m_{2}}) \cdots (1 - r_{t}^{m_{t}}) \cdots > 0$, at any given stage $t$ the Controller has at least $m_{t+1}$ actions available. $\Box$
\end{exl}

The second sufficient condition, unlike the first one, places no restrictions on the goal, and allows the primitive distribution to change across time. Instead, it requires that the number of actions at all stages be dominated by a single random variable with finite mean. Recall that $\#p_{t}$ is the distribution of the cardinality of the action set under the primitive distribution $p_{t}$.

\begin{thm}\label{thm.omni:dominance}
Suppose that there exists a distribution $q$ on $\{1,2,\ldots\}$ with finite mean that dominates $\#p_{t}$ for every $t \in \N$. Then $v^\infty(p_{0},p_{1},\ldots,G) \in \{0,1\}$.
\end{thm}

\begin{thm}\label{thm.all:dominance}
Suppose that there exists a distribution $q$ on $\{1,2,\ldots\}$ with finite mean that dominates $\#p_{t}$ for every $t \in \N$. Then $v^{m}(p_{0},p_{1},\ldots,G) \in \{0,1\}$ for each $m \in \N$. 
\end{thm}

We remark that a time-invariant primitive distribution $p_{0} = p_{1} = \cdots$ satisfies the condition of Theorems \ref{thm.omni:dominance} and \ref{thm.all:dominance} only when $\#p_{0}$ has finite mean, a requirement that Theorems \ref{thm.omni:shift} and \ref{thm.finite:shift} do not impose. 


The two sets of sufficient conditions above imply the zero-one law for both the omniscient values and all the $m$-foresight values. Theorem \ref{thm.finite:dominance} focuses on a given foresight $m$. The sufficient condition of the theorem relies on the notion of a maximal branching process (MBP, Lamperti \cite{Lamperti70, Lamperti72}). Given a distribution $q$ on $\N$, the MBP generated by $q$ is a Markov process $Y_{0} = 1, Y_{1}, Y_{2}, \ldots$ with transitions defined by the formula
\[\mathbb{P}(Y_{t+1} \leq n|Y_{t}) = (q(1) + \cdots + q(n))^{Y_{t}}.\] 
The process is said to be recurrent if for each $n \in \N$ with $q(n) > 0$, almost surely there are infinitely many $t \in \N$ for which $Y_{t} = n$. 

\begin{thm}\label{thm.finite:dominance}
Let $m \in \N$. Suppose that there exists a distribution $q$ on $\N$ such that (a) the {\rm MBP} generated by $q$ is recurrent, and (b) the distribution $q$ dominates $\#p_{t}^{m}$ for each $t \in \N$. Then $v^{m}(p_{0},p_{1},\ldots,G) \in \{0,1\}$. 
\end{thm}

Whether a given MBP is recurrent is, generally, a non-trivial question. Fortunately, Lamperti \cite[Theorem 1]{Lamperti70} derives a sufficient condition for the recurrence of an MBP in terms of the distribution $q$ generating it: namely 
\begin{equation}\label{eqn.Lamperti_cond}
\limsup_{n \to \infty}n(1 - q(1) -  \cdots - q(n)) < e^{-\gamma},
\end{equation}
where $\gamma$ is Euler's constant. We call a distribution $q$ satisfying \eqref{eqn.Lamperti_cond} a \textit{Lamperti distribution}. Importantly, a distribution with finite mean is a Lamperti distribution.

We proceed with two remarks highlighting the relations between various conditions. 

\begin{rem}\rm\label{rem.conditions0}
It is not hard to see that the distribution $\#p_{t}^{m+1}$ dominates $\#p_{t}^{m}$. Consequently, the requirement of Theorem \ref{thm.finite:dominance} becomes progressively stronger as $m$ increases. It is trivially satisfied if $m = 0$ (for $\#p_{t}^{0}$ is a Dirac measure on 1), so that Proposition \ref{thm.0foresight} follows at once. 
\end{rem} 

\begin{rem}\rm\label{rem.conditions1}
The sufficient condition in Theorem \ref{thm.finite:dominance} is implied by the  sufficient condition of Theorems \ref{thm.omni:dominance} and \ref{thm.all:dominance}. To see this, let $q$ be as in Theorem \ref{thm.all:dominance}, and let $g$ be the corresponding probability generating function. Let $g^{m}$ be the $m$-fold composition $g \circ \cdots \circ g$, and $q^{m}$ the corresponding distribution. Since $q$ has finite mean, so does $q^{m}$, and hence it meets the condition of Theorem \ref{thm.finite:dominance}. In particular, Theorem \ref{thm.all:dominance} is a corollary of Theorem \ref{thm.finite:dominance}.  
\end{rem}

While the condition of Theorem \ref{thm.omni:dominance} implies the condition of Theorem \ref{thm.finite:dominance} for any $m \in \N$, the two conditions are not equivalent, the point our final counterexample illustrates. In this counterexample, there is a distribution $q$ satisfying \eqref{eqn.Lamperti_cond} that dominates $\#p_{t}$ for each $t \in \N$. Thus, the counterexample fulfils the sufficient condition of Theorem \ref{thm.finite:dominance} with $m = 1$. And yet, the omniscient value is between $0$ and $1$. This illustrates that the sufficient condition of Theorem \ref{thm.finite:dominance} is weaker than that of Theorem \ref{thm.omni:dominance}. 

\begin{exl}\label{exl.lyons}\rm
For $t \in \N$ let $p_t$ be the primitive distribution that assigns a positive probability to exactly 13 different elements of $\A$, as follows
\begin{center}
\renewcommand{\arraystretch}{1.4}
\begin{tabular}{lc}\toprule
set & probability\\\midrule
$\{0\}$ & $1 - \tfrac{1}{4}2^{-t} - \ldots - \tfrac{1}{4}2^{-t-11}$\\
$\{0,1,\ldots,2^{t}\}$ & $\tfrac{1}{4}2^{-t}$\\
\hspace{1cm}$\vdots$ & $\vdots$\\
$\{0,1,\ldots,2^{t+11}\}$ & $\tfrac{1}{4}2^{-t-11}$\\\bottomrule
\end{tabular}
\end{center}
Let the goal be $G = \{x \in \N^\N: x_t = 0\text{ for at most finitely many }t \in \N\}$. As shown in the following section, the example satisfies the condition of Theorem \ref{thm.finite:dominance} for $m = 1$; yet $0 < v^\infty(p_{0},p_{1},\ldots,G) < 1$.
\end{exl}

Proofs of the claims in the examples are collected in Section 5. In Section $6$, we prove the results on the omniscient value (Theorems \ref{thm.finite:shift} and \ref{thm.omni:dominance}), and Section $7$ is devoted to the proof of the results on the $m$-foresight value (Theorems \ref{thm.finite:shift}, \ref{thm.all:dominance} and \ref{thm.finite:dominance}, and Proposition \ref{thm.0foresight}).

\section{Proof of the claims in the examples}
\subsection{Example \ref{exl.1foresight:invariant}}
\noindent\textsc{Claim 1:} $v^\infty(p_{0},p_{1},\ldots,G) < 1$. For each $t \in \N$ consider the set $\Lambda = \cap_{t \in \N} \Lambda_{t}$, where $\Lambda_{t} = \{\omega \in \Omega: x_{t} \in  M_{0} \cup \cdots \cup M_{t}\text{ for each }x \in [\omega]\}$. As $\Lambda$ is disjoint from $\Omega_{G}$, it suffices to prove that it has a positive probability under the measure $\mu_{0}$. We argue that
\[\mu_{0}(\Lambda) \geq \prod_{t \in \N}r_{t}^{\prod_{i < t}m_{i}}.\]
The desired conclusion then follows from the condition (b) of the example.

Recalling that $\omega_{t}$ denotes the set of nodes of $\omega$ of length $t$, we note that
\[\mu_{0}(\Lambda_{t}|\mathscr{F}_{t}) = \mu_{0}(\{A_{\omega}(h) 
\in \{M_{0},\ldots,M_{t}\}\text{ for every }h \in \omega_{t}\} |\mathscr{F}_{t}) = r_{t}^{\#\omega_{t}}.\]
Let $t \geq 1$. We know that $\#\omega_{t} \leq m_{0} \cdot \cdots \cdot m_{t-1}$ for each $\omega \in \Lambda_{0} \cap \cdots \cap \Lambda_{t-1}$, and consequently $\mu_{0}(\Lambda_{t}|\mathscr{F}_{t}) \geq r_{t}^{m_{0} \cdot \cdots \cdot m_{t-1}}$ everywhere on $\Lambda_{0} \cap \cdots \cap \Lambda_{t-1}$. 
But this implies that 
\[\mu_{0}(\Lambda_{t}|\Lambda_{0} \cap \cdots \cap \Lambda_{t-1}) \geq r_{t}^{m_{0} \cdot \cdots \cdot m_{t-1}},\]
from which the desired estimate follows.\smallskip

\noindent\textsc{Claim 2:} $0 < v^{1}(p_{0},p_{1},\ldots,G)$. Consider the following strategy $\sigma$: at a stage $n \in \N$, choose the action with the largest number of follow-up actions (break the ties by picking the smallest action). Consider the goal $E = \{x \in \N^\N: x_{t} \in M_{t+1}\cup M_{t+2} \cup \cdots\text{ for each }t \in \N\}$, a subset of $G$. We argue that $\sigma$ succeeds for the goal $E$ with positive probability. We accomplish this by showing that
\[\mu_{0}(\Omega_{E}^{\sigma}) \geq \prod_{t \in \N}(1 - r_{t}^{m_{t}}).\]
The  conclusion then follows by condition (c) of the example.

Let $h_{t}(\omega) \in \N^t$ be the node of the decision tree $\omega$ reached by $\sigma$ at stage $t$ (so in particular $h_{0}(\omega) = \o$). Let $Y_{t}(\omega)$ denote the number of actions available to the Controller at stage $t$, thus $Y_{t}(\omega) = \# A_{\omega}(h_{t}(\omega))$. Note that the random variable $Y_{t} : \Omega \to \N$ takes values in the set $\{m_{0},m_{1},\ldots\}$, and that $Y_{t}(\omega) = m_{n}$ precisely when the action set at node $h_{t}(\omega)$ is $M_{n}$ (the numbers $m_{n}$ are all distinct by condition (a)). Let $\Theta_{t} = \{\omega \in \Omega: Y_{t}(\omega) \leq m_{t}\}$. The event $\cap_{t \in \N}\Theta_{t}^{c}$ is exactly $\Omega_{E}^{\sigma}$, the event that $\sigma$ succeeds for the goal $E$.  

For $t = 0$ we have $\mu_{0}(\Theta_{0}^{c}) = 1 - r_{0} = 1 - r_{0}^{m_{0}}$. For $t \geq 1$ we have:
\[\mu_{0}(\Theta_{t} |\mathscr{F}_{t}) = r_{t}^{Y_{t-1}}.\]
To see this, note that the event $\Theta_{t}$ occurs exactly when none of the actions at node $h_{t-1}(\omega)$ could be followed up by more than $m_{t}$ actions. On $\Theta_{t-1}^{c}$, we have $m_{t-1} < Y_{t-1}$, and consequently $m_{t} \leq Y_{t-1}$. We therefore obtain
\[\mu_{0}(\Theta_{t} |\mathscr{F}_{t}) \leq r_{t}^{m_{t}}\text{ everywhere on }\Theta_{t-1}^{c}.\]
But this means that 
\[\mu_{0}(\Theta_{t} | \Theta_{t-1}^{c} \cap \cdots \cap \Theta_{0}^{c}) \leq r_{t}^{m_{t}}.\]
We conclude that 
\[\mu_{0}(\Omega_{E}^{\sigma}) = \prod_{t \in \N} \mu_{0}(\Theta_{t}^{c} | \Theta_{t-1}^{c} \cap \cdots \cap \Theta_{0}^{c}) \geq  \prod_{t \in \N}(1 - r_{t}^{m_{t}}).\]

\subsection{Example \ref{exl.lyons}}
\noindent\textsc{Claim 1:} $v^\infty(p_{0},p_{1},\ldots,G) < 1$. The proof is identical to that in Example \ref{exl.1foresight:shift}.\medskip

\noindent\textsc{Claim 2:} $0 < v^\infty(p_{0},p_{1},\ldots,G)$. We prove a stronger statement, namely that $0 < \linebreak v(p_{0},p_{1},\ldots,G')$, where the goal $G'$ is given by $G' = \{x \in \N^\N: x_t \neq 0\text{ for each }t \in \N\}$. 

Given $\omega \in \Omega$ let $z_{t}(\omega)$ be the cardinality of the set $\omega \cap \{1,2,\ldots\}^{t}$, i.e. the number of nodes of $\omega$ in generation $t$ consisting of actions with no $0$s. As is easy to see, $\{z_{t}\}_{t \in \N}$ is a BPVE,  its offspring in generation $t \in \N$ is as follows: $0$ children with probability $1-\tfrac{1}{4}\sum_{i = 0}^{11}2^{-t-i}$, and, for each $i \in \{0, \ldots, 11\}$, $2^{t+i}$ children with probability $\tfrac{1}{4}2^{-t-i}$.

Clearly, if $z_{t}(\omega) > 0$ for each $t \in \N$ then (by K\"{o}nig's lemma) $[\omega]$ has an infinite sequence of positive natural numbers, so $ \omega \in \Omega_{G'}$. Thus, it suffices to show that the process $\{z_{t}\}_{t \in \N}$ does not go extinct almost surely. We apply the classical sufficient condition due to Fearn \cite[Corollary 4]{Fearn71}. The offspring distribution of $z_{t}$ has mean $\bar{m}_{t} = 3$ and variance \[\bar{\sigma}_{t}^2 \leq \sum_{i = 0}^{11} \tfrac{1}{4}2^{-t-i} \cdot 2^{2t+2i} \leq \tfrac{1}{4}\cdot12\cdot2^{11}\cdot2^{t}.\] 
And therefore
\[\sum_{t \in \N}\frac{\bar{\sigma}_{t}^2}{\bar{m}_{1} \cdots \bar{m}_{t-1} \bar{m}_{t}^2} \leq \tfrac{1}{4}\cdot12\cdot2^{11}\cdot \sum_{t \in \N}\frac{2^{t}}{3^{t}} < \infty,\]
establishing the result.\medskip

\noindent\textsc{Claim 3}: The example satisfies the condition of Theorem \ref{thm.finite:dominance} for $m = 1$ (and consequently its $1$-foresight value is $0$). To see this, let $q$ be the distribution on $\{1,2,\ldots\}$ given as follows: $q(1) = \tfrac{1}{2}$ and $q(1+2^{t}) = \tfrac{1}{4}2^{-t}$ for each $t \in \N$. Then $q$ dominates $\#p_{t} = \#p_{t}^{1}$ for each $t \in \N$ and satisfies Lamperti's condition \eqref{eqn.Lamperti_cond}, since
\[\limsup_{n \to \infty}n(1 - q(1) -  \cdots - q(n)) \leq \limsup_{t \to \infty} (1 + 2^{t})\sum_{i \in \N}\tfrac{1}{4}2^{-t-i} = \tfrac{1}{2} < e^{-\gamma}.\]
Hence it fulfils the requirements (a) and (b) of Theorem \ref{thm.finite:dominance} with $m = 1$. $\Box$

\section{Proof of the results on the omniscient value}
We prove Theorems \ref{thm.omni:shift} and \ref{thm.omni:dominance}. Suppose that, contrary to the conclusion of the theorems, $0 < \mu_{0}(\Omega_{G}) < 1$. We use the following implication of the Levy's zero-one law (see Appendix):
\begin{align}
\mu_{0}(\Omega_{G}) &= \mu_{0}(\{\mu_{0}(\Omega_{G}|\mathscr{F}_{t}) \to 1\})\label{eqn.Levy0},\\ 
\mu_{0}(\Omega_{G}^c) &= \mu_{0}(\{\mu_{0}(\Omega_{G}|\mathscr{F}_{t}) \to 0\}).
\label{eqn.Levy1}\end{align}

Using the fact that $G$ is a tail set and property (2) of $\mu_{0}$ we obtain 
\begin{align}
\mu_{0}(\Omega_{G}^c|\mathscr{F}_{t}) &= \mu_{0}(\{\omega \in \Omega :[\omega] \cap G = \o\} |\mathscr{F}_{t})\nonumber\\ &= \mu_{0}(\{\omega \in \Omega: [\omega(h)] \cap G_{t} = \o\text{ for each }h \in \omega_{t}\}|\mathscr{F}_{t})\nonumber\\ &= \mu_{0}(\{\omega \in \Omega: \omega(h) \in \Omega_{G_{t}}^c\text{ for each }h \in \omega_{t}\}|\mathscr{F}_{t})\nonumber\\ &= (\mu_{t}(\Omega_{G_{t}}^c))^{\#\omega_{t}}\nonumber\\ &= s_{t}^{\#\omega_{t}},\label{eqn.cond_prob}
\end{align}
where we define
\[s_{t} = \mu_t(\Omega_{G_{t}}^c).\]

Taking the expectation in \eqref{eqn.cond_prob} with respect to the measure $\mu_{0}$, we get \[s_0 = \mu_{0}(\Omega_{G}^c) = \mathbb{E}[s_{t}^{\#\omega_{t}}] = g_{0}^{t}(s_t),\]
where the last equality uses the fact that $g_{0}^{t}$ is the probability generating function of $\#\omega_{t}$ under $\mu_{0}$. Since $g_{0}^{t+1} = g_{0}^{t} \circ g_{t}$, it follows that
\begin{equation}\label{eqn.recur}
g_{t}(s_{t+1}) = s_{t}
\end{equation}
for each $t \in \N$.

In particular, the sequence $s_{0} \leq s_{1} \leq \cdots$ is non-decreasing. Moreover, we must have
\begin{equation}\label{eqn.st->1}
\sup_{t \in \N}s_{t} = 1.
\end{equation}
To see this, define $s_{*} = \sup_{t \in \N}s_{t}$ and suppose that $s_{*} < 1$. Since $1 \leq \#\omega_{t}$ for each $t \in \N$ and each $\omega \in \Omega$, equation \eqref{eqn.cond_prob} implies that $\mu_{0}(\Omega_{G}^c|\mathscr{F}_{t}) \leq s_{*} < 1$ on $\Omega$. Equation \eqref{eqn.Levy1} implies that $\mu(\Omega_{G}^{c}) = 0$, contradicting our supposition.\medskip

\noindent\textsc{Proof of Theorem \ref{thm.omni:shift}}: We let $k_{0} = 0$ (the general case can be handled similarly) and write $k$ for $k_{1}$. Then $G_{t} \subseteq G_{t+k}$ and so $s_{t+k} \leq s_{t}$ for each $t \in \N$. We conclude that $s_{0} = s_{1} = s_{2} = \cdots$, in contradiction to \eqref{eqn.st->1}.\medskip

\noindent\textsc{Proof of Theorem \ref{thm.omni:dominance}}: Suppose first that for each $t \in \N$ the primitive distribution $p_{t}$ only assigns positive probability to singleton action sets. Since then $\#\omega_{t} = 1$ almost $\mu_{0}$-surely, we have $\mu_{0}(\Omega_{G}^c|\mathscr{F}_{t}) = s_{t}$, which means that for a given $t \in \N$ the function $\mu_{0}(\Omega_{G}^c|\mathscr{F}_{t})$ is constant on $\Omega$. Invoking \eqref{eqn.Levy0} and \eqref{eqn.Levy1} we obtained the desired contradiction.

Henceforth suppose that for some $t \in \N$ the primitive distribution $p_{t}$ assigns a positive probability to some non-singleton action set. Without loss of generality we may assume that $p_{0}$ fulfills the requirement.

Let $q$ be the distribution as in the hypothesis of the theorem, $m$ its mean, and $g$ its probability generating function. We have $\lim_{x \uparrow 1}g'(x) = m$. Since $g \leq g_{t}$, \eqref{eqn.recur} implies that $g(s_{t+1}) \leq s_{t} \leq s_{t+1}$ for each $t \in \N$. Using these facts and \eqref{eqn.st->1}, we obtain 
\begin{equation}\label{enq.ratiologs}
1 \leq \lim_{t \to \infty}\frac{\ln(s_{t})}{\ln(s_{t+1})} \leq \lim_{t \to \infty}\frac{\ln(g(s_{t+1}))}{\ln(s_{t+1})} = \lim_{t \to \infty}\frac{s_{t+1}g'(s_{t+1})}{g(s_{t+1})} = m.
\end{equation}

Let $\Theta$ denote the event $\{\omega \in \Omega: s_{t}^{\#\omega_{t}} \to 1\}$. We have
\begin{align*}
\Theta &= \{\omega \in \Omega: \#\omega_{t}\ln(s_{t}) \to 0\}\\ &= \{\omega \in \Omega: \#\omega_{t-1}(a)\ln(s_{t}) \to 0\text{ for each }a \in \omega_{1}\}\\ &= \{\omega \in \Omega: \#\omega_{t}(a)\ln(s_{t}) \to 0\text{ for each }a \in \omega_{1}\}\\ &= \{\omega \in \Omega: \omega(a) \in \Theta\text{ for each }a \in \omega_{1}\},
\end{align*}
where we write $\omega_t(a)$ to denote $(\omega(a))_t$, the $t$th generation of the subtree of $\omega$ following action $a$ at stage $0$, and likewise for $\omega_{t-1}(a)$. Here the second equality stems from the fact that 
\[\#\omega_{t}(a) = \sum_{a \in \omega_{1}}\#\omega_{t-1}(a),\] 
and the third is justified by \eqref{enq.ratiologs}. 
Therefore,
\[\mu_{0}(\Omega_{G}^{c}|\mathscr{F}_{1}) = \mu_{0}(\Theta|\mathscr{F}_{1}) =  \mu_{0}(\Theta)^{\#\omega_{1}} =  \mu_{0}(\Omega_{G}^{c})^{\#\omega_{1}}.\]
Here the first and the third equations hold since, by Levy's zero-one law, $\Theta$ coincides with $\Omega_{G}^{c}$ up to a set of $\mu_{0}$-measure zero, and the middle equality is by property (2) of the measure $\mu_{0}$.

Taking the expectations (with respect to $\mu_{0}$), we obtain $s_{0} = g_{0}(s_{0})$. By our supposition (that $p_{0}$ assigns a positive probability to some non-singleton action set), the function $g_{0}$ only has $0$ and $1$ as fixed points. Hence $s_{0}$ is $0$ or $1$, yielding a contradiction to the initial supposition.

\section{Proof of the results on the $m$-foresight value}
\subsection{Main ideas of the proof}\label{subsecn.outline}
This subsection, written deliberately in a non-rigorous manner, aims to convey the gist of the arguments. The three main ingredients of the proof are the notion of a value of a CP at a particular state, the notion of the size of the state, and the bounds on the probability of transition from one state to another.

In Subsection \ref{subs.markov}, we cast the scenario of $m$-foresight as a Markov decision process (MDP). This allows us to make use of the notion of the value of the MDP at a particular state. The states of the MDP in question represent the information that is revealed to the Controller at a particular stage. Recall that, under $m$-foresight, at each stage, Nature reveals to the Controller the first $m+1$ generations of the subtree of the decision tree reached thus far. We think of this information as a \textit{state}. Accordingly, the set of states of the MDP is defined as $S_{m} = S = \{\omega_{m+1}:\omega \in \Omega\}$. Following the Controller's action $a_{t}$ at state $s_{t} \in S$, the transition probabilities to a new state $s_{t+1}$ reflect two considerations: the first $m$ generations of $s_{t+1}$ are completely determined by $s_{t}$ and $a_{t}$, while the $(m+1)$st generation of $s_{t+1}$ is constructed using the action sets that are chosen, independently, from the primitive distribution $p_{t+m+1}$. This gives rise to an MDP with a countable state space, finite action sets, and non-homogenous transitions.

Given a state $s \in S$ we define $v_{t}^m(s)$ to be the value of the MDP, conditional on $s \in S$ being the state revealed to the Controller at stage $t$. 

Let $\phi_{t}$ denote the measure on $S$ given by $\phi_{t}(Q) = \mu_{t}(\{\omega \in \Omega: \omega_{m+1} \in Q\})$\label{nota.phit} for each $Q \subseteq S$. This is the \textit{unconditional distribution of the state $s_{t}$} at stage $t \in \N$. Consider some $\epsilon> 0$ and a stage $t$. We partition the set of states into two subsets, those having stage $t$ value smaller than $\epsilon$ (called \textit{bad}), and those with the stage $t$ value at least $\epsilon$ (called \textit{good}). The set of good states at stage $t$ is thus $Z_{t}^{\epsilon} = \{s \in S: \epsilon \leq v_{t}^{m}(s)\}$\label{nota.Zt}. Let $z_{t}^{\epsilon} = \phi_{t}(Z_{t}^{\epsilon})$\label{nota.zt} be the unconditional probability for a state at stage $t$ to be good.

Theorem \ref{thm.finite:shift} and Theorem \ref{thm.finite:dominance} follow from the three propositions stated below and proven in Subsection \ref{subs.dichotomy}. Theorem \ref{thm.all:dominance} and Proposition \ref{thm.0foresight} are both implied by Theorem \ref{thm.finite:dominance}.


\begin{prop}\label{pro.dichotomyI}
Suppose that there exists an $\epsilon > 0$ such that $\limsup_{t \to \infty} z_{t}^{\epsilon} > 0$. Then the $m$-foresight value of the {\rm CP} is $1$.
\end{prop}

\begin{prop}\label{pro.dichotomyII}
Suppose that the primitive distribution is time-invariant and that there exist $0 \leq k_0 < k_1$ such that $G_{k_{0}} \subseteq G_{k_{1}}$. If $\limsup_{t \to \infty} z_{t}^{\epsilon}=0$ for every $\epsilon>0$, then the $m$-foresight value of the {\rm CP} is $0$.
\end{prop}

\begin{prop}\label{pro.dichotomyIII}
Suppose that there exists a distribution $q$ on $\N$ such that (a) the {\rm MBP} generated by $q$ is recurrent, and (b) the distribution $q$ dominates $\#p_{t}^{m}$ for each $t \in \N$. If $\limsup_{t \to \infty} z_{t}^{\epsilon}=0$ for every $\epsilon>0$, then the $m$-foresight value of the {\rm CP} is $0$.
\end{prop}

In a nutshell, the key to the proof of Proposition \ref{pro.dichotomyI} is to look at the sets of states that the Controller can guarantee to reach in $m+1$ stages; here we are particularly interested in a transition from the state $s_{t}$ at stage $t$ to a good state at stage $t+m+1$, i.e. to a state in the set $Z_{t+m+1}^{\epsilon}$. Proposition \ref{pro.dichotomyII} poses no major difficulty. The key to Proposition \ref{pro.dichotomyIII} is to look at the conditional probability that, given a state $s_{t}$ at stage $t$, a state at stage $t+m+1$ is bad. We seek to bound this probability from below, uniformly over all $m$-foresight strategies. This leads us to introduce the notion of the size of the state.\medskip 

\noindent\textsc{Idea of the proof of Proposition \ref{pro.dichotomyI}:} Fix $\epsilon  > 0$ witnessing the hypothesis. The key argument of the proof is that starting at a stage $t$, the Controller is able to guarantee a transition to the set of states $Z_{t+m+1}^{\epsilon}$ by stage $t+m+1$ with probability $z_{t+m+1}^{\epsilon}$. In fact, this the Controller can do using any $0$-foresight strategy. It follows that the value is bounded away from zero uniformly across the states. This is only possible if all states have value 1.\medskip

\noindent\textsc{Idea of the proof of Proposition \ref{pro.dichotomyII}:}  Suppose first that, as in Theorem \ref{thm.finite:dominance}, the primitive distribution is time-invariant, and also that - for the sake of clarity - the goal is shift-invariant. In this case $Z_{t}^{\epsilon}$ and $z_{t}^{\epsilon}$ are independent of $t$, and thus $z_{0}^{\epsilon} = 0$ for each $\epsilon \in (0,1)$, and hence, almost surely, the value of a state at stage $0$ has value $0$. The argument is essentially the same under the hypothesis of the proposition.\medskip

\noindent\textsc{Idea of the proof of Proposition \ref{pro.dichotomyIII}:} Let us illustrate the main ideas in a very special setup. Suppose that the primitive distribution is time-invariant, that C has 1-foresight (so $m = 1$), and assume that she adopts the 1-maximizing strategy (recall Example \ref{exl.benchmark}).

Under these suppositions, the number of actions available to the Controller at a particular stage forms a MBP induced by the primitive distribution. Under the hypothesis of the Proposition \ref{pro.dichotomyIII}, this process is recurrent. And thus, almost surely, the number of actions available to the Controller has a bounded subsequence. 

Now consider the probability of a 2-step transition to a bad state, that is, the probability of transition from a state $s_{t}$ at stage $t$ to a state in $S \setminus Z_{t+2}^{\epsilon}$ at stage $t+2$. The exact transition probability is not easy to compute, because the action at stage $t+2$ is not measurable with respect to $s_{t}$. Instead, we rely on a lower bound, obtained as follows: if \textit{each} action $a_{t+1}$ available to the Controller at stage $t+1$ leads to a subtree that induces a bad state, a visit to a bad state is clearly inevitable. The probability for any given subtree at stage $t+2$ to induce a bad state is $1-z_{t+2}^{\epsilon}$. By independence, these probabilities can be multiplied. We thus obtain a lower bound on the probability of transition to a bad state, depending on (a) the unconditional probability for a state to be bad, and (b) the number of actions available to the Controller at stage $t+1$.

Since, under the hypothesis of the Proposition \ref{pro.dichotomyIII}, the unconditional probability of a good state $z_{t}^{\epsilon}$ approaches $0$, while the number of available actions is bounded, the probability of transition to a bad state is bounded away from zero. And hence, the Controller almost surely encounters infinitely many bad states. But this implies that her probability of success, under the strategy she employs, is zero. 

The proof of Proposition \ref{pro.dichotomyIII} generalizes these ideas in several ways.

\medskip   

The rest of the section is structured as follows. The next subsection provides the details on the Markov decision process and a recaps a few known facts about the function $v_{t}$. The proof of Propositions \ref{pro.dichotomyI}--\ref{pro.dichotomyIII} is carried out in Subsection \ref{subs.dichotomy} and is split into 8 steps. 

Step \ref{step.trans:0} treats the $m+1$ step transitions under $0$-foresight strategies. Step \ref{step.v==1} establishes Proposition \ref{pro.dichotomyI}, and Step 3 proves Proposition \ref{pro.dichotomyII}. Steps 4-8 complete the proof of Proposition \ref{pro.dichotomyIII}. 

Step \ref{step.involuntary:trans} derives a lower bound on the probability of transition from stage $t$ to a given set of states $Q \subseteq S$ by stage $t+m+1$.  Step \ref{step.value:upperbound} applies the result of Step \ref{step.involuntary:trans} to the set of bad states $S \setminus Z_{t+m+1}^\epsilon$ to obtain an upper bound on the value at a state $s_{t}$. Step \ref{step.xi:distribution} shows that the distribution of the size of a state is dominated by the distribution of the MBP induced by $q$. Step \ref{step.coupling} provides the details on the coupling of the two processes, and deduces that there is a subsequence of states of bounded size. The final step puts it all together and completes the proof.

\subsection{The $m$-foresight as an MDP}\label{subs.markov}
\noindent\textbf{The MDP:} One can think of the scenario of $m$-foresight as a non-homogenuous Markov decision process (MDP). This view is helpful as we can rely on a number of well known properties of the value.

Define the Markov decision process MDP$^{m}$\label{nota.MDP} as having the following ingredients:  
\begin{itemize}
\item \textit{The state space:} The (countable) state space is $S = S_{m}$\label{nota.S}, the set of non-empty finite subsets of $\N^{m+1}$; 
\item \textit{Action sets:} The (finite) action set $A(s)$\label{nota.A(s)} at state $s \in S$ is the set of natural numbers $a_0 \in \N$ such that $(a_0,a_1,\ldots,a_m) \in s$ for some sequence $(a_1,\ldots,a_m) \in \N^{m}$;
\item \textit{Initial distribution:} The initial state, i.e. the state $s_{0}$ at stage $0$, is distributed according to the measure $\phi_{0}$. 
\item \textit{Transition kernels:} Given a state $s_{0} \in S$ and an action $a_{0} \in A(s_{0})$, transitions with positive probability occur only to states $s_{1} \in S$ with the following property: if $(a_1,\ldots,a_{m+1}) \in s_{1}$, then $(a_{0},a_{1},\ldots,a_{m}) \in s_{0}$. The probability of transition at stage $t$ from $(s_{0},a_{0})$ to such a state $s_{1}$ is 
\[\prod_{\substack{(a_{1},\ldots,a_{m}) \in \N^m:\\(a_{0},a_{1},\ldots,a_{m}) \in s_{0}}}p_{t+m+1}(\{a_{m+1} \in \N: (a_{1},\ldots,a_{m},a_{m+1}) \in s_{1}\}).\]
\item \textit{Payoffs: }The payoff is $1$ if the sequence of actions $(a_{0},a_{1},\ldots)$ is in $G$, and is $0$ otherwise.
\end{itemize}

A strategy in MDP$^{m}$ is a sequence $\sigma^{*} = (\sigma_{0}^{*}, \sigma_{1}^{*}, \ldots)$ where $\sigma_{t}^{*} : S^{t+1} \to \N$. Each strategy in MDP$^{m}$ induces a strategy $\sigma \in \Sigma_{m}$ with $m$-foresight in the original CP, and vice versa. Henceforth we identify a strategy in MDP$^{m}$ with the corresponding $m$-foresight strategy of the CP, and write $\sigma$ for both.  

Let $\Omega^{*} = (S \times \N)^\N$\label{nota.Omega*}. Define $G^{*} = \{(s_{0},x_{0},s_{1},x_{1},\ldots) \in \Omega^{*}: (x_{0},x_{1},\ldots) \in G\}$\label{nota.G*}. For $t \in \N$ we let $\I_{t}$\label{nota.It} denote the sigma-algebra on $\Omega^{*}$ generated by the projection map $\Omega^{*} \to S \times (S \times \N)^{t}$, $(s_{0},x_{0},s_{1},x_{1},\ldots) \mapsto (s_{0},x_{0},\ldots,x_{t-1},s_{t})$.  With a slight abuse of notation we write $s_{t}$ to denote the projection map $\Omega^{*} \to S$, as well as an element of $S$. Thus in particular $\I_{0}$ is the sigma-algebra on $\Omega^{*}$ generated by $s_{0}$. 

A strategy $\sigma \in \Sigma_{m}$ together with the initial distribution and the transition kernels induce a probability measure, denoted $\mu_{\sigma}^{*}$\label{nota.mu*}, on the Borel sigma-algebra of the product space $\Omega^{*}$. Let $\mathbb{E}_{\sigma}^{*}$ be the expectation with respect to the measure $\mu_{\sigma}^{*}$.\medskip

\noindent\textbf{The value of the MDP:} Define the \textit{stage $t$ value}\label{nota.vtm} of the MDP$^{m}$ as the function $v_{t}^{m} : \Omega^{*} \to [0,1]$ given by 
\begin{equation}\label{eqn.markov:value}
v_{t}^{m} = \sup_{\sigma \in \Sigma_{m}}\mu_{\sigma}^{*}(G^{*}|\I_{t}).
\end{equation}
We discuss some properties of the value of the MDP$^{m}$.  

Firstly, as the set $G^{*}$ is tail, the stage $t$ value $v_{t}^{m}$ of MDP$^{m}$ is measurable with respect to $s_{t}$, the state at stage $t$. Henceforth, we write $v^m_{t}(s_{t})$ to denote the stage $t$ value of MDP$^{m}$ at a particular state $s_{t} \in S$.

Secondly, for each $\sigma \in \Sigma$ the process $(v^m_{t})_{t\in\N}$ is a supermartingale on the measure space $(\Omega^{*}, (\mathscr{I}_{t})_{t \in \N}, \mu_{\sigma}^{*})$. That is, for each $t \in \N$
\begin{equation}\label{eqn.martingale}
\mathbb{E}_{\sigma}^{*}(v^m_{t+1}|\mathscr{I}_{t}) \leq v^m_{t}.
\end{equation}
In particular, it converges $\mu_{\sigma}^{*}$-almost surely to the function $v^m_{\infty} : \Omega^{*} \to[0,1]$ given by $v^m_{\infty} = \limsup_{t \to \infty} v^m_{t}$, and 
\begin{equation}\label{eqn.martingale0}
\mathbb{E}_{\sigma}^{*}(v^m_{\infty}|\I_{0}) \leq v^m_{0}.
\end{equation}

The supermartingale property of the value does not rely on the particular structure of MDP$^{m}$ and is well-known in the literature (see for example Dubins, Savage, Sudderth, Gilat \cite[Theorem 1]{DubinsSavageSudderthGilat14} who use the term  “excessive”).

Thirdly, it holds that 
\begin{equation}\label{eqn.value2}
v^m_{0} \leq \sup_{\sigma \in \Sigma_{m}} \mu_{\sigma}^{*}(\{v^m_{\infty} = 1\} | \I_{0}).
\end{equation}
To see \eqref{eqn.value2}, take any $\sigma \in \Sigma_{m}$. By L\'{e}vy's zero-one law $\mu_{\sigma}^{*}(G^{*}|\I_{t}) \to 1_{G^{*}}$ almost surely with respect to the measure $\mu_{\sigma}^{*}$. In view of \eqref{eqn.markov:value}, whenever $\mu_{\sigma}^{*}(G^{*}|\I_{t})$ converges to 1, also $v^m_{t}$ converges to 1, and hence $v^m_{\infty} = 1$. And thus $G^{*} = \{\mu_{\sigma}^{*}(G^{*}|\I_{t}) \to1\} \subseteq \{v^m_{\infty} = 1\}$. Hence $\mu_{\sigma}^{*}(G^{*}|\I_{0}) \leq \mu_{\sigma}^{*}(\{v^m_{\infty} = 1\}|\I_{0})$. Taking the supremum over $\sigma \in \Sigma_{m}$ yields \eqref{eqn.value2}.\medskip

\noindent\textbf{The connection between $\Omega^{*}$ and $\Omega$:} Let $\sigma : \Omega \to \N^{\N}$ be an $m$-foresight strategy. Given a stage $n \in \N$, and a decision tree $\omega \in \Omega$, let $s_{n}^{\sigma}(\omega)$ be the $(m+1)$st generation of the subtree of $\omega$ starting at the node $h_{n}^{\sigma}(\omega)$, thus $s_{n}^{\sigma}(\omega) = (\omega(h_{n}^{\sigma}(\omega)))_{m+1}$\label{nota.stsigma}. This defines the function $s_{n}^{\sigma} : \Omega \to S$. Thus $s_{n}^{\sigma}$ is the state the  Controller encounters at stage $n$, provided that she adheres to the strategy $\sigma$.

Now consider the map $e_{\sigma} : \Omega \to \Omega^{*}$, $\omega \mapsto (s_{0}^{\sigma}(\omega), \sigma_{0}(\omega), s_{1}^{\sigma}(\omega), \sigma_{1}(\omega), \ldots)$. The measure $\mu_{\sigma}^{*}$ on $\Omega^{*}$ is the pullback of the $\mu_{0}$ on $\Omega$ under the map $e_{\sigma}$: for each Borel set $B \subseteq \Omega^{*}$ it holds that $\mu_{\sigma}^{*}(B) = \mu_{0}(e_{\sigma}^{-1}(B))$. Note also that $e_{\sigma} : (\Omega,\F_{t+m+1}) \to (\Omega^{*},\I_{t})$ is measurable.         

\subsection{The proof of Propositions \ref{pro.dichotomyI}, \ref{pro.dichotomyII}, \ref{pro.dichotomyIII}}\label{subs.dichotomy}
Fix the depth of foresight $m \in \N$.

Step \ref{step.trans:0} studies $(m+1)$-step transition probabilities in MDP$^{m}$ under $0$-foresight strategies. Its key insight is that under a $0$-foresight strategy, the state $s_{t+m+1}$ at stage $t+m+1$ is independent of the information available at stage $t$, $\I_{t}$; in particular, the states $s_{t+m+1}$ and $s_{t}$   are independent. To see this, notice that the actions prescribed by the $0$-foresight strategy at stages $t,\ldots,t+m$ are measurable with respect to $s_{t}$, which is determined by the first $t+m+1$ generations of the decision tree. The state $s_{t+m+1}$, on the other hand, is determined by the subsequent $m+1$ generations, i.e. $t+m+2, \ldots, t+2m+2$.\smallskip 

\begin{step}{ (transitions under $0$-foresight strategies)}\label{step.trans:0}
\textit{For any $0$-foresight strategy $\sigma$, each set $Q \subseteq S$ of states, and each stage $t \in \N$ it holds that}
\[\mu_{\sigma}^{*}(\{s_{t+m+1} \in Q\}|\I_{t}) = \phi_{t+m+1}(Q).\]

It suffices to show that
\[\mu_{0}(\{s_{t}^{\sigma} \in Q\}|\mathscr{F}_{t}) = \phi_{t}(Q)\]
for each $t \in \N$. The desired equality then follows: indeed, the measure $\mu_{\sigma}^{*}$ is the pullback of $\mu_{0}$ under the map $e_{\sigma}$. The preimage of the set $\{s_{t+m+1} \in Q\} \subseteq \Omega^{*}$ under $e_{\sigma}$ is the set $\{s_{t+m+1}^{\sigma} \in Q\} \subseteq \Omega$, and the preimage of any element of the collection $\I_{t}$ under $e_{\sigma}$ is an element of $\F_{t+m+1}$.  

Recall that $s_{t}^{\sigma}(\omega)$ is the $m+1$st generation of $\omega(h_{t}^{\sigma}(\omega))$, the subtree of $\omega$ starting at the node $h_{t}^{\sigma}(\omega) = (\sigma_{0}(\omega), \ldots, \sigma_{t-1}(\omega))$. Since $\sigma$ has $0$-foresight, the functions $\sigma_{0}, \ldots, \sigma_{t-1}$ are measurable with respect to $\mathscr{F}_{t}$, and so is $\omega \mapsto h_{t}^{\sigma}(\omega)$. Therefore, under the conditional measure $\mu_{0}(\cdot|\mathscr{F}_{t})$, the distribution of $\omega(h_{t}^{\sigma}(\omega))$ is $\mu_{t}$. And thus $\mu_{0}(\{s_{t}^{\sigma} \in Q\}|\mathscr{F}_{t}) = \mu_{t}(\{\omega_{m+1} \in Q\}) = \phi_{t}(Q)$.
\end{step}

The intuition behind the next step is the following: Applying Step \ref{step.trans:0} to the set $Q = Z_{t+m+1}^{\epsilon}$ of good states, we conclude that the Controller, using any $0$-foresight strategy, can guarantee that a good state is visited with probability $z_{t+m+1}^{\epsilon}$. Hence, if some subsequence of the sequence $z_{t}^{\epsilon}$ is bounded away from zero, she can be sure to be in a good state infinitely often. And this implies that the value of the CP is $1$.\medskip

\begin{step}{ (proof of Proposition \ref{pro.dichotomyI})}\label{step.v==1}
\textit{Suppose that there is an $\epsilon > 0$ such that $\limsup_{t \to \infty}z_{t}^{\epsilon} > 0$. Then the $m$-foresight value of the CP is $1$.}

Since $z_{t}^{\epsilon}$ is non-increasing in $\epsilon$ for each $t \in \N$, we can choose a small enough $\epsilon > 0$ so that $\limsup_{t \to \infty}z_{t}^{\epsilon} \geq \epsilon$. For each $\sigma \in \Sigma$,
\begin{align*}
v^m_{t}(s_{t}) &\geq \mathbb{E}_{\sigma}^{*}(v^m_{t+m+1}(s_{t+m+1})|\mathscr{I}_{t})\\ &\geq \epsilon \cdot \mu_{\sigma}^{*}(\{s_{t+m+1} \in Z_{t+m+1}^{\epsilon}\}|\mathscr{I}_{t}) \\&= \epsilon \cdot \phi_{t+m+1}(Z_{t+m+1}^{\epsilon})\\&= \epsilon \cdot z_{t+m+1}^{\epsilon},     
\end{align*}
where the inequality in the line is by the supermartingale property of the value \eqref{eqn.martingale}, the inequality of the second line follows by the definition of the set $Z_{t+m+1}^{\epsilon}$, the equality of the third line is by Step \ref{step.trans:0}, and the final equality is by definition of $z_{t+m+1}^{\epsilon}$. It follows that $v^m_{\infty} \geq \epsilon^2$ everywhere on $\Omega^*$. And hence
\begin{align}\label{eqn.maj}
v^m_\infty & \geq \epsilon^2 1_{\{v^m_\infty<1\}}+ 1_{\{v^m_\infty=1\}} = \epsilon^2 + (1-\epsilon^2)1_{\{v^m_\infty=1\}} 
\end{align}

By \eqref{eqn.martingale0}, \eqref{eqn.maj}, and \eqref{eqn.value2}, we obtain
\begin{align*}
v^m_{0}(s_{0}) &\geq \sup_{\sigma \in \Sigma} \mathbb{E}_{\sigma}^{*}(v^m_{\infty}|\I_{0})\\ &\geq\epsilon^2 + (1-\epsilon^2) \cdot \sup_{\sigma \in \Sigma}\mu_{\sigma}^{*}(\{v^m_{\infty} = 1\} | \I_{0})\\ &\geq \epsilon^2 + (1-\epsilon^2) \cdot v^m_{0}(s_{0}),
\end{align*}
which implies that $v^m_{0}(s_{0}) = 1$. Thus $v^m_{0}$ is identically $1$ on $S$, and hence the $m$-foresight value of the CP is $1$ for any $m$.
\end{step}

Proposition \ref{pro.dichotomyII} now follows easily: as the goal periodically becomes more and more permissive, the probability of a good state has a non-decreasing subsequence. And hence, the assumption of the proposition implies that good states have probability zero, implying that the value is zero.\smallskip

\begin{step}{ (proof of Proposition \ref{pro.dichotomyII})} 
\textit{Under the hypothesis of Theorem \ref{thm.finite:shift} if $\limsup_{t \to \infty} z_{t}^{\epsilon} = 0$ for every $\epsilon > 0$, then the $m$-foresight value of the {\rm CP} is $0$.}

As $G_{k_{0}} \subseteq G_{k_{1}}$, it holds that $G_{k_{0} + n} \subseteq G_{k_{1} + n}$ for each $n \in \N$. It can be seen directly from the definition of the value \eqref{eqn.markov:value} that $v_{k_{0} + n}(s) \leq v_{k_{1} + n}(s)$ for each $s \in S$. Therefore $z_{k_{0} + n}^{\epsilon} \leq z_{k_{1} + n}^{\epsilon}$ for each $\epsilon > 0$. Hence $z_{k_{0}}^{\epsilon} \leq \limsup_{t \to \infty} z_{t}^{\epsilon}$ and therefore $z_{k_{0}}^{\epsilon} = 0$ for each $\epsilon > 0$. Hence $v_{k_{0}}(s_{k_{0}}) = 0$ for any state $s_{k_{0}}$ that has a positive probability under the distribution $\phi_{k_{0}}$. The result follows.  
\end{step}


The rest of the section is devoted to the proof of Proposition \ref{pro.dichotomyIII}. The key is the notion of the size of a state. In the case of $1$-foresight (i.e. when $m=1$), the size of a state $s_{t}$ is defined as the largest number of stage $t+1$ actions $a_{t+1}$ that could be followed up on the action $a_{t}$ at stage $t$. In general, the size of the state $s_{t}$ is the largest number of sequences  $a_{t+1},\ldots,a_{t+m}$ of actions that the Controller could take at stages $t+1,\ldots,t+m$, following her action at stage $t$. Formally, we define the \textit{size of the state} $s \in S$ by\label{nota.size}
\[u(s) = \max_{a_0 \in \N}\#\{(a_1,\ldots,a_m) \in \N^m: (a_0,a_1,\ldots,a_m) \in s\}.\]

We use the notion of the size of a state to obtain a lower bound on $(m+1)$-step transition probabilities. This is accomplished by Step \ref{step.involuntary:trans}. More precisely, it provides a lower bound on the conditional probability that, given the information $\I_{t}$ available at stage $t$, the state $s_{t+m+1}$ at stage $t+m+1$ is in some set $Q \subseteq S$. The bound is uniform over all the $m$-foresight strategies. It depends on the unconditional probability for a state at stage $t+m+1$ to be in $Q$, and on the  size of the state $s_{t}$. In particular, the smaller the size, the tighter is the bound. 

It is instructive to contrast transitions under $0$-foresight strategies (studied in Step \ref{step.trans:0} above) and under $1$-foresight strategies (e.g. Step \ref{step.involuntary:trans} below). Under a $0$-foresight strategy, both actions $a_{t}$ and $a_{t+1}$ are measurable with respect to the information at stage $t$, hence the states $s_{t}$ and $s_{t+2}$ are independent, and the conditional probability for the state $s_{t+2}$ to be in $Q$ equals the unconditional probability $\phi_{t+2}(Q)$. Under a $1$-foresight strategy, only the action at stage $t$ is, in general, measurable with respect to $\I_{t}$, while the follow-up action $a_{t+1}$ might depend on the information the Controller will receive at the subsequent stage. Consequently, the states $s_{t}$ and $s_{t+2}$ are no longer independent. 

In the case of 1-foresight, we obtain a lower bound on the conditional probability that $s_{t+2} \in Q$ given $\I_{t}$ as follows. Suppose that the state at stage $t$ is $s_{t}$, and that the Controller chooses action $a_{t}$. If each of the available actions $a_{t+1}$ following $a_{t}$ leads to a subtree inducing a state in $Q$ at stage $t+2$, then clearly the Controller is bound to visit the set $Q$. And since the Controller has at most $u(s_{t})$ actions at stage $t+1$, the probability to visit $Q$ at stage $t+2$ is at least 
\[\phi_{t+2}(Q)^{u(s_{t})}.\] 

\begin{step}{ (a lower bound on transitions)}\label{step.involuntary:trans}
\textit{For each strategy $\sigma \in \Sigma_{m}$, stage $t \in \N$, and the set of states $Q \subseteq S$ it holds that 
\[\mu_{\sigma}^{*}(\{s_{t+m+1} \in Q\}|\mathscr{I}_{t}) \geq \phi_{t+m+1}(Q)^{u(s_{t})}.\]
}
By the same logic as in Step \ref{step.trans:0}, it suffices to argue that 
\[\mu_{0}(\{s_{t+m+1}^{\sigma} \in Q\}|\mathscr{F}_{t+m+1}) \geq \phi_{t+m+1}(Q)^{u(s_{t}^{\sigma})}.\]
This inequality is a consequence of three considerations. Firstly, $\{s_{t+m+1}^{\sigma} \in Q\}$ is the event that the $(m+1)$-st generation of the subtree of $\omega$ starting at the node $h_{t+m+1}^{\sigma}(\omega)$ is in $Q$. Consider the event that for each node $h \in \omega_{t+m+1}$ extending $h_{t+1}^{\sigma}(\omega)$, the $(m+1)$-st generation of $\omega(h)$ is in $Q$. Clearly, the latter event is a subset of the former. Secondly, under the measure $\mu_{0}(\cdot|\mathscr{F}_{t+m+1})$, the random variables $\omega(h)$, for $h \in \omega_{t+m+1}$, are independent and distributed according to $\mu_{t+m+1}$. Thirdly, there are at most $u(s_{t}^{\sigma}(\omega))$ nodes $h \in \omega_{t+m+1}$ extending $h_{t+1}^{\sigma}(\omega)$. This is because a sequence $h \in \omega_{t+m+1}$ extends $h_{t+1}^{\sigma}(\omega)$ if and only if $h$ follows up on $h_{t+1}^{\sigma}(\omega)$ with a sequence of actions $(a_{t+1},\ldots,a_{t+m})$ such that $(\sigma_{t}(\omega),a_{t+1},\ldots,a_{t+m})$ is in $s_{t}^{\sigma}(\omega)$.

Putting these considerations together, we obtain:
\begin{align*}
&\mu_{0}(\{\omega \in \Omega: s_{t+m+1}^{\sigma}(\omega) \in Q\} | \mathscr{F}_{t+m+1}) \geq\\&\mu_{0}(\{\omega \in \Omega: (\omega(h))_{m+1} \in Q\text{ for each }h \in \omega_{t+m+1}\text{ extending }h_{t+1}^{\sigma}(\omega)\} | \mathscr{F}_{t+m+1}) \geq\\ & \mu_{t+m+1}(\{\omega \in \Omega: \omega_{m+1} \in Q\})^{u(s_{t}^{\sigma})} =\\ & \phi_{t+m+1}(Q)^{u(s_{t}^{\sigma})}.
\end{align*}
\end{step}

Step \ref{step.value:upperbound} applies the result of Step \ref{step.involuntary:trans} by letting $Q$ be the set of bad states  $S \setminus Z_{t+m+1}^{\epsilon}$. This leads easily to an upper bound on the value at state $t$.\smallskip

\begin{step}{ (an upper bound on the value)}\label{step.value:upperbound} 
\textit{For each state $s \in S$ and stage $t \in \N$
\[v_{t}^{m}(s) \leq 1 - \epsilon \cdot (1-z_{t+m+1}^{\epsilon})^{u(s)}.\]}

Essentially this follows from the fact that, by Step \ref{step.involuntary:trans}, conditional on the state $s$ being revealed to the Controller at stage $t$, the probability to visit the set of states $S \setminus Z_{t+m+1}^{\epsilon}$ at stage $t+m+1$ is at least $(1-z_{t+m+1}^{\epsilon})^{u(s)}$. 

To derive the result formally, take a strategy $\sigma \in \Sigma$. Consider also the subsets $I = \{s_{t} = s\}$ and $J = \{s_{t+m+1} \in S \setminus Z_{t+m+1}^{\epsilon}\}$ of $\Omega^{*}$. We have
\begin{align*}
\mu_{\sigma}^{*}(G^{*}|I) &\leq 1 - \mu_{\sigma}^{*}(J|I) + \mu_{\sigma}^{*}(G^{*} \cap J|I)\\ &= 1 - \mu_{\sigma}^{*}(J|I) + \mu_{\sigma}^{*}(J|I) \cdot \mu_{\sigma}^{*}(G^{*}|J \cap I)\\ &\leq 1 - \mu_{\sigma}^{*}(J|I) + \mu_{\sigma}^{*}(J|I) \cdot \epsilon\\ &= 1 - (1 - \epsilon) \cdot \mu_{\sigma}^{*}(J|I)\\*[-4pt] &\leq 1 - (1 - \epsilon) \cdot (1-z_{t+m+1}^{\epsilon})^{u(s)},
\end{align*}
where the inequality of the third line is by Equation \eqref{eqn.markov:value} and the fact that $v_{t+m+1}(s_{t+m+1}) \leq \epsilon$ everywhere on the event $J$. The final inequality is by Step \ref{step.involuntary:trans}. 

The desired inequality is now obtained by taking the supremum with respect to $\sigma$ and using Equation \eqref{eqn.markov:value}.
\end{step}

The point of the following two steps is to establish that, under any $m$-foresight strategy, the process $u(s_{0}), u(s_{m}), u(s_{2m}), \ldots$ is dominated by an MBP induced by $q$. Since we assume the MBP to be recurrent, this immediately implies that almost surely there is a subsequence of states bounded in size. 

The argument is clearest in the benchmark case of 1-foresight. Consider the event that the size of the state $s_{t+1}$ is at most $n$. This happens exactly when any action $a_{t+1}$ that the Controller might take at stage $t+1$ could be followed up by at most $n$ actions at stage $t+2$. the Controller has at most $u(s_{t})$ actions at stage $t+1$. And for any given of these actions, the probability that the next stage action set has cardinality at most $n$ is $p_{t+2}(1)+ \ldots + p_{t+2}(n)$. Finally, the primitive distribution $p_{t+2}$ is dominated by $q$. All this implies that the probability of the event $u(s_{t+1}) \leq n$ conditional on $\I_{t}$ is at least 
\[(q(1)+ \ldots + q(n))^{u(s_{t})}.\]

\begin{step}{ (the size of a state is dominated by an MBP)}\label{step.xi:distribution}
\textit{For each $t \in \N$,
\[\mu_{\sigma}^{*}(\{u(s_{t+m}) \leq n\} |\mathscr{I}_{t}) \geq (q(1)+\cdots+q(n))^{u(s_{t})}.\]}

By the same logic as in Step \ref{step.trans:0}, it suffices to argue that
\[\mu_{0}(\{u(s_{t+m}^{\sigma}) \leq n\} |\F_{t+m+1}) \geq (q(1)+\cdots+q(n))^{u(s_{t}^{\sigma})}.\]

Firstly, by the definition of $X$ we have  
\[u(s_{t}^{\sigma}(\omega)) = \max_{a_{t} \in \N} \#\{(a_{t+1},\ldots,a_{t+m}) \in \N^{m}: (a_{t},a_{t+1},\ldots,a_{t+m}) \in s_{t}^{\sigma}(\omega)\}.\]
Now recall that $s_{t}^{\sigma}(\omega)$ is the $m+1$st generation of the subtree of $\omega$ starting at the node $h_{t}^{\sigma}(\omega)$. A sequence $(a_{t},a_{t+1},\ldots,a_{t+m})$ is in the $m+1$st generation of the subtree $\omega(h_{t}^{\sigma}(\omega))$ if and only if $(h_{t}^{\sigma}(\omega),a_{t}) \in \omega$ and the sequence $(a_{t+1},\ldots,a_{t+m})$ is in the $m$th generation of the subtree $\omega(h_{t}^{\sigma}(\omega),a_{t})$. We can thus express $u(s_{t}^{\sigma}(\omega))$ as 
\begin{equation}\label{eqn.xi:def}
u(s_{t}^{\sigma}(\omega))  = \max_{a_{t} \in \N:\,(h_{t}^{\sigma}(\omega),a_{t}) \in \omega} \#\omega_{m}(h_{t}^{\sigma}(\omega),a_{t}).
\end{equation}
Since $h_{t+1}^{\sigma}(\omega)$ is a node of $\omega$ extending $h_{t}^{\sigma}(\omega)$ we obtain 
\begin{equation}\label{eqn.xi:boundlower}
\#\omega_{m}(h_{t+1}^{\sigma}(\omega)) \leq u(s_{t}^{\sigma}(\omega)).
\end{equation}

Consider now $u(s_{t+m}^{\sigma}(\omega))$. Each node $(h_{t+m}^{\sigma}(\omega),a_{t+m})$ of $\omega$ is a concatenation of the node $h_{t+1}^{\sigma}(\omega)$ with some element $h$ of $\omega_{m}(h_{t+1}^{\sigma}(\omega))$, the $m$th generation of the subtree $\omega(h_{t+1}^{\sigma}(\omega))$. And thus from \eqref{eqn.xi:def} there follows the inequality
\begin{equation}\label{eqn.xi:boundupper}
u(s_{t+m}^{\sigma}(\omega)) \leq \max_{h \in \omega_{m}(h_{t+1}^{\sigma}(\omega))} 
\#\omega_{m}(h_{t+1}^{\sigma}(\omega),h).
\end{equation}

We obtain: 
\begin{align*}
\mu_{0}(\{u(s_{t+m}^\sigma) \leq n\} |\mathscr{F}_{t+m+1}) \geq& \mu_{0}\left(\left\{\begin{array}{c} \#\omega_{m}(h_{t+1}^{\sigma}(\omega),h) \leq n\\ \text{ for each }h \in \omega_{m}(h_{t+1}^{\sigma}(\omega)) \end{array}\right\} |\mathscr{F}_{t+m+1} \right)\\
=&\mu_{t+m+1}(\{\#\omega_{m} \leq n\})^{\#\omega_{m}(h_{t+1}^{\sigma})}\\ \geq& \mu_{t+m+1}(\{\#\omega_{m} \leq n\})^{u(s_{t}^\sigma)}\\ \geq& (q(1)+\cdots+q(n))^{u(s_{t}^\sigma)}.
\end{align*}
Here the inequality in the first line is by \eqref{eqn.xi:boundupper}. The equation of the second line holds because $h_{t+1}^{\sigma}(\omega)$ is measurable with respect to $\mathscr{F}_{t+m+1}$, the $m$th generation of $\omega(h_{t+1}^{\sigma}(\omega))$ is also measurable with respect to $\mathscr{F}_{t+m+1}$, and because under the measure $\mu_{0}(\cdot|\mathscr{F}_{t+1})$ the random variables $\omega(h')$, $h' \in \omega_{t+m+1}$, are independent and distributed according to $\mu_{t+m+1}$. The inequality of the third line is by \eqref{eqn.xi:boundlower}. The last inequality is by the hypothesis of the theorem.
\end{step}

\begin{step}{ (there is a subsequence of states of bounded size)}\label{step.coupling}
\textit{For each $\sigma \in \Sigma$,}
\[\textstyle\liminf_{t \to \infty}u(s_{t}) < \infty,\, \mu_{\sigma}^{*}\textit{-almost surely}.\]


The key is the inequality of Step \ref{step.xi:distribution}. Using coupling (more precisely, Lindvall \cite[Theorem 8.5]{Lindvall02}), one can show that there exist random sequences $(U_{k}')_{k \in \N}$ and $(Y_{k}')_{k \in \N}$ such that $U_{k}' \leq Y_{k}'$ for each $k \in \N$, the first sequence is distributed like the stochastic process $(u(s_{mk}))_{k \in \N}$ under the measure $\mu_{\sigma}^{*}$, and the second one like the MBP $(Y_{k})_{k \in \N}$ generated by the distribution $q$. Since the maintained hypothesis the MBP $(Y_{k})_{k \in \N}$ generated by $q$ is recurrent, the result follows. 

To apply the said theorem, one sets up two transition kernels from $\N^{k-1}$ to $\N$, as follows. Consider the probability of transition from a point of $\N^{k-1}$, say $(n_{0},\ldots,n_{k-1})$, to a point of $\N$, for example $n_{k}$. Under the first transition kernel, the corresponding transition occurs with probability 
\[\mu_{\sigma}^{*}(\{u(s_{mk}) = n_{k}\} \mid \{u(s_{0}) = n_{0},\ldots,u(s_{mk-m}) = n_{k-1}\}).\] Under the second one, it occurs with probability $\mathbb{P}(Y_{k} = n_{k}|Y_{k-1} = n_{k-1})$, or, more explicitly, $(q(1) + \cdots + q(n_{k}))^{n_{k-1}} - (q(1) + \cdots + q(n_{k}-1))^{n_{k-1}}$. 

The inequality of Step \ref{step.xi:distribution}, applied with $t = mk-m$ implies that the second transition kernel dominates the first.
\end{step}

We are in a position to complete the proof. Since there is, almost surely, a subsequence of steps of bounded size, and since the unconditional probability of a good state approaches zero, the Controller will visit a bad state infinitely often. This implies that the value of the MDP is zero.\medskip 

\begin{step}{ (proof of Proposition \ref{pro.dichotomyIII})} 
Let $\epsilon = \tfrac{1}{2}$ (we can take any $\epsilon \in (0,1)$). We know that $z_{t}^{\epsilon} \leq \tfrac{1}{2}$ for large $t \in \N$. Take any $\sigma \in \Sigma_{m}$. By Step \ref{step.value:upperbound} then $v_{t}(s) \leq 1 - 2^{-u(s)}$. Hence the event$\{v_{\infty} = 1\}$ of $\Omega^{*}$ is contained in the event $\{\liminf_{t \to \infty}u(s_{t}) = \infty\}$. But by Step \ref{step.coupling}, the latter has probability $0$ under $\mu_{\sigma}^{*}$, and hence so does the former. Inequality \eqref{eqn.martingale0} now implies that $v_{0}(s_{0}) = 0$ almost $\mu_{\sigma}^{*}$-surely. Hence the CP has value $0$.
\end{step}

\section{Concluding remarks}
\subsection{Counterexamples and the tightness of results}
\begin{table}
\renewcommand{\arraystretch}{1.4}
\caption{Summary of the counterexamples}\label{table.counterexamples}
\begin{center}
\begin{tabular}{lccc}\toprule
& Example \ref{exl.1foresight:shift} & Example \ref{exl.1foresight:invariant} & Example \ref{exl.lyons}\\\midrule
omniscient value & in $(0,1)$ & in $(0,1)$  & in $(0,1)$\\
1-foresight value & in $(0,1)$ & in $(0,1)$ & $0$\\
goal shift-invariant? & Yes & No & Yes\\
primitive distribution time-invariant? & No & Yes & No\\
the mean of primitive distributions & & &\\*[-5pt]
uniformly bounded? & No & No & Yes\\
primitive distributions dominated by & & &\\*[-5pt]
a Lamperti distribution? & No & No & Yes\\\bottomrule
\end{tabular}
\end{center}
\end{table}

Table \ref{table.counterexamples} gives the summary of our three counterexamples.

Examples \ref{exl.1foresight:shift} and \ref{exl.1foresight:invariant} illustrate the tightness of the conditions of Theorems \ref{thm.omni:shift} and \ref{thm.finite:shift}. They witness that neither shift-invariance of goal, nor time-invariance of the primitive distribution, on their own, suffice for the zero-one laws. 

Example \ref{exl.lyons} illustrates the discrepancy between the sufficient conditions of Theorems \ref{thm.omni:dominance} and \ref{thm.finite:dominance}. In the example, the omniscient value is strictly between 0 and 1, even though the mean of the distribution $\#p_{t}$ equals $3$ for each $t \in \N$. In fact, the distributions $p_{t}$ are finitely supported, and so all the moments of $\#p_{t}$ are finite. As we have argued above, there does exist a Lamperti distribution (i.e. the distribution satisfying the inequality \eqref{eqn.Lamperti_cond}) dominating $\#p_{t}$ for each $t \in \N$,  thus fulfilling the condition of Theorem \ref{thm.finite:dominance}. However, any such distribution necessarily has infinite mean.

\subsection{One-foresight and MBPs}
This subsection is intended as a response to Lamperti's \cite{Lamperti70} invitation for the readers (see the third paragraph of the introduction) to suggest applications of MBPs. The scenario of 1-foresight is intimately related to the notion of MBP, the link that Theorem \ref{thm.finite:dominance} exploits. To further highlight the connection, we focus on the situation when primitive distributions are time-invariant ($p_{t} = p$ for all $t \in \N$), and suppose that the Controller adopts 1-step maximizing strategy $\sigma$ as in Example \ref{exl.benchmark}. Under $\sigma$, the number of actions available to the Controller at stage $t+1$ equals the size of the state at stage $t+1$: is $\#A_{\omega}(h_{t}^{\sigma}(\omega)) = u(s_{t}^{\sigma}(\omega))$. And both these processes are MBPs induced by the distribution $p$. 

Note that we have made no claim of optimality of the 1-step maximizing strategy. It would be interesting to know in what classes of CPs this strategy is actually optimal.   
\subsection{The omniscient value and rates of growth of the decision tree}
We have seen two sufficient conditions (Theorems \ref{thm.omni:shift} and \ref{thm.omni:dominance}) for the omniscient value to obey the zero-one law. Below we state yet another condition in terms of the rate of growth of the size of the decision tree.

\begin{lemma}\label{thm.omni:rate}
Let $r_{0}, r_{1},\ldots$ be a sequence of positive numbers such that the process $r_{t}^{-1} \cdot \#\omega_{t}$ converges to a positive finite limit $\mu_{0}$-almost surely. Then $v^{\infty}(p_{0},p_{1},\ldots,G) \in \{0,1\}$.
\end{lemma}
\begin{proof}
By \eqref{eqn.Levy0}--\eqref{eqn.Levy1}, the process $\mu_{0}(\Omega_{G}^c|\mathscr{F}_{t})$ converges to $0$ with probability $\mu_{0}(\Omega_{G})$ and to $1$ with probability $\mu_{0}(\Omega_{G}^c)$. In particular, no point other than $0$ and $1$ is an accumulation point of the process with positive probability. By \eqref{eqn.cond_prob} we have
\[\mu_{0}(\Omega_{G}^c|\mathscr{F}_{t}) = \left(s_{t}^{r_{t}}\right)^{r_{t}^{-1} \cdot \#\omega_{t}}.\]
Since $r_{t}^{-1} \cdot \#\omega_{t}$ converges to a positive finite limit, it must be the case that $s_{t}^{r_{t}} \to 0$ or $s_{t}^{r_{t}} \to 1$. In the first case $\mu_{0}(\Omega_{G}) = 1$, and in the second case $\mu_{0}(\Omega_{G}) = 0$.
\end{proof}

Recall that the process $\{\#\omega_{t}\}_{t \in \N}$ is a BPVE, with $\{\#p_{t}\}_{t \in \N}$ the corresponding sequence of offspring distributions. The sequence $\{r_{t}\}_{t \in \N}$ appearing in the statement of Lemma \ref{thm.omni:rate} is called a \textit{normalizing sequence}. Whether a BPVE admits a normalizing sequence is one of the central questions in the theory of branching processes.

By the classical result of Heyde \cite{Heyde70}, the normalizing sequence as required by Lemma \ref{thm.omni:rate} does exist if the distributions $\#p_{t}$ do not depend on time and have finite mean. By Lyons \cite[Theorem 4.14]{Lyons92}, it is sufficient that the support of the distributions $\{\#p_{t}\}_{t \in \N}$ be uniformly bounded. Both of these conditions are subsumed by the hypothesis of Theorem \ref{thm.omni:dominance}. 

Another sufficient condition is Kersting's condition (A), regularity, applied to the sequence of distributions $\#p_0, \#p_1, \ldots$ (see Kersting \cite[Theorem 2(ii)]{Kersting20}). Kersting's regularity is logically unrelated to the condition of Theorem \ref{thm.omni:dominance}. Yet another sufficient condition for the existence of the normalizing sequence has been developed in d'Souza and Biggins \cite[Theorem 2]{dSouza92}. We do not know whether there is a relationship between this latter condition and the condition of Theorems \ref{thm.omni:dominance}.

We also do now know whether the sufficient condition of Theorem  \ref{thm.omni:dominance} (that there exists a single distribution with finite mean dominating $\#p_{t}$ for each $t \in \N$) is sufficient for the existence of the normalizing sequence.  

\section{Appendix}
\noindent\textbf{$\Omega_{G}$ is analytic:} We argue that $\Omega_{G}$ is an analytic subset of $\Omega$.

Recall that, identifying a decision tree $\omega$ with its body $[\omega]$, we identify $\Omega$ with $\K = \K(\N^{\N})$, the space of non-empty compact subsets of $\N^{\N}$ with the Vietoris topology. Since under the said identification $\Omega_{G}$ is mapped onto $\K_{G} = \{K \in \K: K \cap G \neq \o\}$, it suffices to show that the latter is an analytic subset of $\K$. 

To see this, consider the set $\J_{G} = \{(x,K) \in \N^\N \times \K: x \in K \cap G\}$. The set $\J_{G}$ is a Borel subset of $\N^\N \times \K$, as it is the intersection of the closed set $\{(x,K) \in \N^\N \times \K: x \in K\}$ (Kechris \cite[Exercise 4.29]{Kechris95}) and the Borel set $G \times \K$. Finally, $\K_{G}$ is the projection of $\J_{G}$ onto the second coordinate.\medskip

\noindent\textbf{The decision-theoretic definition of the omniscient value:}
For each strategy $\sigma$ it holds that $\Omega_{G}^{\sigma} \subseteq \Omega_{G}$. Conversely, there exists a strategy $\sigma$ such that $\Omega_{G}^{\sigma} = \Omega_{G}$. 

Indeed, by Jankov--von Neumann uniformization theorem (Kechris \cite[Theorem 18.1]{Kechris95}) the set $\J_{G}$ admits an analytically-measurable uniformizing function, i.e. a function $\sigma : \K_{G} \to \N^\N$ such that $(K, \sigma(K)) \in \J_{G}$ for each $K \in \K_{G}$. Now extend $\sigma$ to all of the space $\K$ by letting $\sigma(K)$ be the leftmost branch of $K \in \K$. The map $\K \to \N^{\N}$ mapping the set $K$ to its leftmost branch is continuous. Hence $\sigma$ is the strategy with the desired property.\medskip

\noindent\textbf{Levy's zero-one law:} As $\Omega_{G}$ is analytic, it differs from a Borel subset of $\Omega_{G}$, say $\Theta \subseteq \Omega_{G}$, by a set of $\mu_{0}$-measure zero (Bertsekas and Shreve \cite[Lemma 7.26]{BertsekasShreve96}). Hence $\mu_{0}(\Omega_{G}|\mathscr{F}_{t}) = \mu_{0}(\Theta|\mathscr{F}_{t}) \to \mu_{0}(\Theta|\mathscr{F}) = 1_{\Theta} = 1_{\Omega_{G}}$, almost surely with respect to the measure $\mu_{0}$ (the convergence is by e.g. Bogachev \cite[Example 10.3.14]{Bogachev07}). The equalities \eqref{eqn.Levy0} and \eqref{eqn.Levy1} follow.

\end{document}